\documentclass[10pt]{article}

\usepackage{amsmath}
\usepackage{amsfonts}
\usepackage{amssymb}
\usepackage{amsthm}
\usepackage[T1]{fontenc}
\usepackage[latin1]{inputenc}
\usepackage[greek,francais,english]{babel}
\usepackage{graphics}
\usepackage{vmargin}
\setmarginsrb{1in}{1in}{1in}{1in}{0cm}{0cm}{0cm}{1.5cm}

\newcommand{\R}{\mathbb R}
\newcommand{\N}{\mathbb N}
\newcommand{\Z}{\mathbb Z}
\newcommand{\K}{\mathbb K}
\newcommand{\e}{\varepsilon}
\newcommand{\be}{\begin{equation}}
\newcommand{\ee}{\end{equation}}

\newcommand{\p}{\partial}

\newcommand{\fe}{f^{\varepsilon}}
\newcommand{\geps}{g^{\varepsilon}}
\newcommand{\re}{r^{\varepsilon}}

\newcommand{\sgn}{\mathrm{sgn}}

\newcommand{\supp}{\mathrm{Supp}}
\newcommand{\Id}{\mathrm{Id}}

\newcommand{\mean}[1]{\left\langle #1\right\rangle}

\newcommand{\om}{\omega}
\newcommand{\xs}{\xi^{\sharp}}

\newtheorem{prop}{Proposition}[section]
\newtheorem{lem}{Lemma}[section]
\newtheorem{thm}{Theorem}
\newtheorem{defi}{Definition}[section]
\newtheorem{corol}{Corollary}[section]
\newtheorem{rem}{Remark}[section]
\begin{document}
\title{Homogenization of linear transport equations in a stationary ergodic setting}

\author{Anne-Laure Dalibard\footnote{CEREMADE-UMR 7534,
Université Paris-Dauphine, Place du maréchal de Lattre de
Tassigny, 75775 Paris Cedex 16, FRANCE; e-mail address :
\texttt{dalibard@ceremade.dauphine.fr}}}
\date{}
\bibliographystyle{amsplain}
\maketitle

\begin{abstract}
We study the homogenization of a linear kinetic equation which
models the evolution of the density of charged particles submitted
to a highly oscillating electric field. The electric field and the
initial density are assumed to be random and stationary. We
identify the asymptotic microscopic and macroscopic profiles of
the density, and we derive formulas for these profiles when the
space dimension is equal to one.

\end{abstract}

\section{Introduction}

This note is concerned with the homogenization of a linear
transport equation in a stationary ergodic setting. The equation
studied here describes the evolution of the density of charged
particles in a rapidly oscillating random electric potential. This
equation can be derived by passing to the semi-classical limit in
the Schrödinger equation (see \cite{Gerard_semiclass},
\cite{LionsPaul}, and the presentation in \cite{HamdacheFrenod}).
Our work generalizes a result of E. Frénod and K. Hamdache (see
\cite{HamdacheFrenod}) which was obtained in a periodic setting.
The strategy of proof we have chosen here is different from the
one of \cite{HamdacheFrenod}, and allows us to retrieve some of
the results in \cite{HamdacheFrenod} in a rather simple and
explicit fashion.

Let us mention a few related works on the homogenization of linear
transport equations; we emphasize that this list is by no means
exhaustive. In \cite{AmiratHamdacheZiani}, Y. Amirat, K. Hamdache
and A. Ziani study the homogenization of a linear transport
equation in a periodic setting and give an application to a model
describing a multidimensional miscible flow in a porous media. In
\cite{DumasGolse} (see also \cite{Golsetransport}), Laurent Dumas
and François Golse focus on the homogenization of linear transport
equations with absorption and scattering terms, in periodic and
stationary ergodic settings. And in \cite{wetransport}, Weinan E
derives strong convergence results for the homogenization of
linear and nonlinear transport equations with oscillatory
incompressible velocity fields in a periodic setting.

Let us now present the context we will be working in : let
$(\Omega, \mathcal F,P)$ be a probability space, and let
$(\tau_x)_{x\in\R^N}$ be a group transformation acting on
$\Omega$. We assume that $\tau_x$ preserves the probability
measure $P$ for all $x\in\R^N$, and the group transformation is
ergodic, which means
$$
\forall A\in\mathcal F,\quad \left(\tau_xA =A\ \forall
x\in\R^N\Rightarrow P(A)=0 \ \text{or}\  1\right).
$$
The periodic setting can be embedded the stationary ergodic
setting (see \cite{PapanicolaouVaradhan}). We will denote by
$E[\cdot]$ the expectation with respect to the probability measure
$P$; in the periodic case, we will write $\mean{f}$ rather than
$E[f]$ to refer to the average of $f$ over one period.

We consider a potential function $u=u(y,\om)\in
L^{\infty}(\R^N\times\Omega)$ which is assumed to be stationary,
i.e.
$$
u(y+z,\om)=u(y,\tau_z\om)\quad \forall
(y,z,\om)\in\R^N\times\R^N\times\Omega
$$
Moreover, we assume that $0\leq u(y,\om)\leq u_{\text{max}}=\sup
u$ for all $y\in\R^N,\om\in\Omega$, and $u(\cdot,\om)\in
W^{2,\infty}_{\text{loc}}(\R^N)$ for almost every $\om\in\Omega$,
so that $\nabla_y u\left(y,\omega \right)$ is well-defined and
locally Lipschitz continuous with respect to its $y$ variable.

Let $\fe=\fe(t,x,\xi,\om)$, ($t\geq 0,\ x\in\R^N, \xi\in\R^N,
\omega\in\Omega$) be the solution of the transport equation
\be\left\{
\begin{array}{l} \p_t \fe(t,x,\xi, \omega) + \xi \cdot \nabla_x
\fe(t,x,\xi, \omega) - \frac{1}{\e} \nabla_y
u\left(\frac{x}{\e},\omega \right)\cdot
\nabla_{\xi}\fe(t,x,\xi,\omega)=0,\\
\fe(t=0,x,\xi,\omega)=f_0\left(x,\frac{x}{\e},\xi,\omega \right).
\end{array}
\right. \label{eq:kinhomog} \ee

Here, we assume that the initial data $f_0=f_0(x,y,\xi,\om)$
belongs to $L^1_{\text{loc}}(\R^N_{x}\times\R^N_{\xi},
L^{\infty}(\R^N_y\times\Omega))$ and is stationary in $y$, i.e.
$$
f_0(x,y+z,\xi,\om)=f_0(x,y,\xi,\tau_z\om)\quad \text{for all }
(x,y,z,\xi,\om)\in\R^{4N}\times \Omega.
$$

It is well-known from the classical theory of linear transport
equations that for every $\om\in\Omega$, there exists a unique
solution $\fe$ of \eqref{eq:kinhomog} in
$L^{\infty}_{\text{loc}}((0,\infty),
L^1_{\text{loc}}(\R^N_x\times\R^N_{\xi}))$. The goal of this paper
is to study the asymptotic behavior of $\fe$ as $\e\to 0$. Thus,
following \cite{HamdacheFrenod}, we define the constraint space
$\mathbb K$ :

\begin{defi}
Let \be \xi\cdot \nabla_y f(y,\xi,\om)- \nabla_y u(y,\om)\cdot
\nabla_{\xi} f(y,\xi,\om)=0 \label{eq:constraint}\ee be the
constraint equation, and let
$$
\mathbb K:=\{f\in L^1_{\text{loc}}(\R^N_{\xi}\times\R^N_y,L^1(
\Omega)); \ f\ \text{satisfies \eqref{eq:constraint} in}\ \mathcal
D'(\R^N_{y}\times\R^N_{\xi})\ \text{a.s. in}\ \om\}.
$$

We also define the projection $P$ onto the constraint space
$\mathbb K$, characterised by $P(f)\in \mathbb K$ for $f\in
L^1_{\text{loc}}(\R^N_{\xi}\times\R^N_y,L^1( \Omega))$ stationary,
and
$$
\int_{\R^N\times\Omega}(P(f)-f)(y,\xi,\om)\:
g(y,\xi,\om)\:d\xi\:dP(\om)=0\quad \text{for a.e.}\ y\in\R^N
$$
for all stationary functions $g\in
L^\infty(\R^N_y\times\R^N_{\xi}\times \Omega)\cap \mathbb K$, with
compact support in $\xi$.

(A more precise definition of the projection $P$ will be given in
the second section).

Finally, we define $\mathbb K^{\bot}$ as
$$
\mathbb K^{\bot}:=\{f\in
L^1_{\text{loc}}(\R^N_{\xi}\times\R^N_y,L^1( \Omega)); \exists g
\in L^1_{\text{loc}}(\R^N_{\xi}\times\R^N_y,L^1( \Omega)),\quad
f=P(g)-g\}.
$$

\end{defi}

\begin{rem}
Let us indicate that the constraint equation can easily be derived
thanks to a formal two-scale Ansatz : indeed, assume that
$$
\fe(t,x,\xi,\om)\approx f\left(t,x,\frac{x}{\e},\xi,\om
\right)\quad \text{as} \ \e\to 0;
$$
inserting this asymptotic expansion in equation
\eqref{eq:kinhomog}, we see that $f$ necessarily satisfies the
constraint equation \eqref{eq:constraint}.
\end{rem}

\begin{rem} Let $f,g\in L^{\infty}(\R^N_y,
L^2(\R^N_{\xi}\times\Omega))$ be stationary, and assume that
$f\in\K$ and $g\in\K^{\bot}$. Then for a.e. $y\in\R^N$,
$$
\int_{\R^N\times\Omega} f(y,\xi,\om)
g(y,\xi,\om)\;d\xi\;dP(\om)=0.
$$
This is a characterization of $\K^{\bot}$ for the class of
stationary functions in $L^{\infty}(\R^N_y,
L^2(\R^N_{\xi}\times\Omega))$.

\end{rem}

Here, we provide another proof for the result of E. Frénod and K.
Hamdache in \cite{HamdacheFrenod} in the ``non-perturbed case''.
Our proof is based on the use of the ergodic theorem, and gives a
more concrete insight of the projection $P$ and of the microscopic
behavior of the sequence $\fe$. Moreover, it allows us to retrieve
the explicit formulas of the integrable case.

The first result we prove in this paper is the following
\begin{thm}
Let $f_0\in L^1_{\text{loc}}(\R^N_x\times \R^N_{\xi}\times\R^N_y;
L^1( \Omega))$ stationary.

Let $\fe=\fe(t,x,\xi,\om)$ be the solution of \eqref{eq:kinhomog}.
Then for all $\e> 0$, there exist $f=f(t,x,y,\xi,\om)$ and
$g=g(t,x;\tau,y,\xi,\om)$, both stationary in $y$, and a sequence
$\{\re(t,x,\xi,\om)\}_{\e>0}$ such that

$$\fe(t,x,\xi,\om)=f\left(t,x,\frac{x}{\e},\xi,\om \right)+ g \left(t,x;\frac{t}{\e},\frac{x}{\e},\xi,\om \right)+ \re(t,x,\xi,\om)$$ and :
\begin{itemize}
\item
$||\re||_{L^1_{\text{loc}}((0,\infty)\times\R^N_x\times\R^N_{\xi},
L^1(\Omega))}\to 0$ as $\e\to 0$;

\item $f\in
L^{\infty}_{\text{loc}}((0,\infty);L^1_{\text{loc}}(\R^N_x\times\R^N_{\xi}\times\R^N_y;
L^1( \Omega)))$, and $f(t,x)\in \mathbb K$ for a.e. $t\geq 0$,
$x\in\R^N$;

\item For all $T>0$, for all compact
$K\subset\R^N_x\times\R^N_\xi\times\R^N_y$,
$$\sup_{0\leq t\leq T,\\0\leq \tau\leq T}||g||_{L^1(K\times \Omega)
}<\infty.$$ Moreover, $g(t,x;\tau,\cdot)\in \mathbb K^{\bot}$ for
a.e. $(t,x,\tau)\in(0,\infty)\times\R^N\times(0,\infty)$;

\item Microscopic evolution equation for $g$ : for a.e.
$t,x\in(0,\infty)\times\R^N$, $g(t,x;\cdot)$ is a solution of \be
\frac{\p g}{\p \tau} + \xi\cdot \nabla_y g -   \nabla_y u \cdot
\nabla_{\xi} g =0.\label{eq:evol_micro}\ee

Moreover, for all $T>0$
$$
\left|\left| \int_0^T g\left(t,x;\frac{t}{\e},\frac{x}{\e},\xi,\om
\right)\:dt\right|\right|_{L^1_{\text{loc}}(\R^N_x\times\R^N_{\xi},
L^1(\Omega))}\to 0 \quad \text{as}\  \e\to 0.
$$

\item Macroscopic evolution equation : $f$ and $g$ satisfy \be
\p_t \left(\begin{array}{c}f\\g\end{array}\right) +
\xi^{\sharp}(y,\xi,\om) \cdot \nabla_x
\left(\begin{array}{c}f\\g\end{array}\right)=0,
\label{eq:evol_homog}\ee where
$$
\xi^{\sharp}(y,\xi,\om):=P(\xi)(y,\xi,\om);
$$

\item Initial data :
\begin{gather*}
f(t=0,x,y,\xi,\om)=P(f_0)(x,y,\xi,\om),\\
g(t=0,x;\tau=0,y,\xi,\om)= \left[f_0-P(f_0)\right](x,y,\xi,\om).
\end{gather*}
\end{itemize}

\label{thm:general}
\end{thm}

Before going any further, we wish to make a few comments on the
above results. First, let us stress that it is not obvious that
the function $g$ is well-defined : indeed, let $S(t)$ ($t\geq 0$)
denote the semi-group associated to the macroscopic evolution
equation \eqref{eq:evol_homog}, and let $T(\tau)$ ($\tau\geq 0$)
be the semi-group associated to the microscopic evolution equation
\eqref{eq:evol_micro}. Then $g$ is well defined if and only if,
for all stationary function $g_0=g_0(x,y,\xi,\om)$, for all
$t,\tau\geq 0$,
$$
T(\tau)\left[S(t) g_0\right]=S(t) \left[T(\tau)g_0 \right].
$$
This identity follows from the fact that the speed
$\xs(y,\xi,\om)$ appearing in equation \eqref{eq:evol_homog} is a
stationary solution of \eqref{eq:evol_micro} by definition of the
projection $P$, and is thus invariant by the semi-group $T(\tau)$.

Next, let us explain briefly the meaning of theorem
\ref{thm:general}. The idea is the following : write $f_0$ as
$f_0=f_{0\parallel} + f_{0\bot}$, with
$f_{0\parallel}(x,\cdot)\in\K$ and $f_{0\bot}(x,\cdot)\in
\K^{\bot}$ a.e. Then $\fe$ can be written as $\fe_{\parallel} +
\fe_{\bot}$, where $\fe_{\parallel}$ (resp. $\fe_{\bot}$) is the
solution of equation \eqref{eq:kinhomog} with initial data
$f_{0\parallel}\left(x,\frac x{\e},\xi,\om \right)$ (resp.
$f_{0\bot}\left(x,\frac x{\e},\xi,\om \right)$). Theorem
\ref{thm:general} states that
$$
\fe_{\parallel}-f\left(t,x,\frac x{\e},\xi,\om \right)\to 0
$$
strongly in $L^1_{\text{loc}}$ norm. In particular, there are no
microscopic oscillations in time in this part of $\fe$. We wish to
emphasize that this result appears to us to be new.

We now focus on the other part, namely $\fe_{\bot}$. An easy
consequence of the theorem is
$$
\int_0^T\fe_{\bot}(t,x,\xi,\om)\;dt \to 0
$$
in $L^{\infty}_{\text{loc}}(\R^N_x; L^1_{\text{loc}}(\R^N_{\xi};
L^1(\Omega)))$ and for all $T>0$. However, it would be wrong to
think that $\fe_{\bot}$ vanishes in
$L^1_{\text{loc}}((0,\infty)\times\R^N_x\times\R^N_{\xi},
L^1(\Omega))$, for instance. Indeed
$$
\fe_{\bot}(t,x,\xi,\om)\approx g\left(t,x;\frac t{\e},\frac
x{\e},\xi,\om \right)
$$
in $L^1_{\text{loc}}$, and
$$
||g(t=0,x;\tau,y)||_{L^1(\R^N_{\xi}\times\Omega)}=||f_{0\bot}(x,y)||_{L^1(\R^N_{\xi}\times\Omega)}
$$
as soon as $f_{0\bot}(x,y)\in L^1(\R^N_{\xi}\times\Omega)$ for
almost every $x,y$. Consequently, if $f_{0\bot}\neq 0$, then for
all $T>0$ and for all compact $K\subset \R^N$, there exists a
constant $C>0$ depending only on $K$,
$||f_{0\bot}||_{L^1_{\text{loc}}}$, and $T$ such that
$$
\int_0^T\left|\left|g\left(t,x;\frac t {\e},\frac x
{\e},\xi,\om\right)\right|\right|_{L^1(K\times\R^N_{\xi}\times\Omega)}\;
dt\geq C.
$$
Hence $\fe_{\bot}$ does not vanish strongly in general. In other
words, there are fast oscillations in time, due to the
ill-preparedness of the initial data (i.e. $f_0(x,\cdot)\notin
\K$), but these oscillations do not cancel out as $\e$ vanishes.

\vskip1mm

Let us now explain briefly here how our strategy of proof differs
from the one of E. Frénod and K. Hamdache. The key of our analysis
lies in the study of the behavior as $\e\to 0$ of the Hamiltonian
system
$$
\left\{
\begin{array}{l}
\dot{Y}^{\e}(t,x,\xi,\om)=-\Xi^{\e}(t,x,\xi,\om),\quad t>0\\
\dot{\Xi}^{\e}(t,y,\xi,\om)=\frac{1}{\e}\nabla_y u(Y^{\e}(t,x,\xi,\om),\om),\quad t>0\\
Y^{\e}(t=0,x,\xi,\om)=x,\ \Xi^{\e}(t=0,x,\xi,\om)=\xi,\quad
(x,\xi,\om)\in\R^{2N}\times\Omega.\end{array} \right.
$$
Indeed,
$$
\fe(t,x,\xi,\om)=f_0\left(Y^{\e}(t,x,\xi,\om),\frac{{Y}^{\e}(t,x,\xi,\om)}{\e},{\Xi}^{\e}(t,y,\xi,\om),\om
\right),
$$
so that we can deduce the asymptotic behavior of $\fe$ from the
one of $({Y}^{\e},\Xi^{\e})$. And it is easily checked that
$$
\begin{array}{l}
{Y}^{\e}(t,x,\xi,\om)=\e Y\left(\frac{t}{\e},\frac{x}{\e},\xi,\om \right)\\
{\Xi}^{\e}(t,y,\xi,\om)=\Xi\left(\frac{t}{\e},\frac{x}{\e},\xi,\om
\right),\end{array}
$$
where $(Y,\Xi)$ is the solution of the system \be \left\{
\begin{array}{l}
\dot{Y}(t,y,\xi,\om)=-\Xi(t,y,\xi,\om),\quad t>0\\
\dot{\Xi}(t,y,\xi,\om)=\nabla_y u(Y(t,y,\xi,\om),\om),\quad t>0\\
Y(t=0,y,\xi,\om)=y,\ \Xi(t=0,y,\xi,\om)=\xi,\quad
(y,\xi,\om)\in\R^{2N}\times\Omega.\end{array} \right.
\label{eq:systhamilt} \ee Hence, in order to study the limit of
$\fe$ as $\e\to 0$, we have to investigate the long time behavior
of the system $(Y,\Xi)$, and this will be achieved with the help
of the ergodic theorem in the second section.

\vskip3mm

In the case when $N=1$, we can give explicit formulas for
$\xs(y,\xi,\om)$; the proof of this formula in the stationary
ergodic case is strongly linked to methods from the Aubry-Mather
theory (see \cite{WEAubryMather}, \cite{EvansGomes1},
\cite{LionsSougCorrectors}), and thus also to the homogenization
of Hamilton-Jacobi equations. Let us first recall the definition
of the homogenized Hamiltonian $\bar{H}$ (see \cite{LPV})
$$
\bar{H}(p)=u_{\text{max}} + \frac 1 2 \left\{\begin{array}{ll}
                                                0\quad&\text{if}\ |p|<
                                                E\left[\sqrt{2(u_{\text{max}}-u)}\right]\\
                                                \lambda\quad&\text{if}\
                                                |p|\geq
                                                E\left[\sqrt{2(u_{\text{max}}-u)}\right],\quad
                                                \text{where }|p|=
                                                E\left[\sqrt{2(u_{\text{max}}-u)+\lambda}\right]\end{array}\right.
$$

\begin{prop}Assume that $N=1$.

Let $(y,\xi,\om)\in\R\times\R\times\Omega$ such that
$H(y,\xi,\om)>u_{\text{max}}$. Assume that for all $Q\in\R$, for
all $\om\in \Omega$, there exists a Lipschitz continuous function
$v(\cdot,\om)$, viscosity solution of
$$
H(y, Q +\nabla_y v(y,\om),\om)=\bar{H}(Q)
$$
such that \be\frac{v(y,\om)}{1 + |y|}\to 0 \quad\text{as}\ |y|\to
\infty \label{hyp:correctsublin} \ee a.s. in $\om$.

Let $P=P(y,\xi,\om)\in \R$ such that $\bar{H}(P)=H(y,\xi,\om)$ and
$sgn(P)=sgn(\xi)$. Then
$$
\xs(y,\xi,\om)=\bar{H}'(P)
$$

Moreover, if $L$ is the dual function of $H$, i.e.
$$
L(y,p,\om)=\sup_{\xi\in\R}\left(p\xi- H(y,\xi\,\om) \right)=\frac
1 2 |p|^2 - u(y,\om),
$$
and $\bar{L}$ is the homogenized Lagrangian, then
$$
P(L)(y,\xi,\om)=\bar{L}(\xs(y,\xi,\om).)
$$
\label{prop:xs_ergo}
\end{prop}

In the periodic case, we will give another proof of the above
result; the strategy chosen in that case is inspired from
techniques and calculations in classical mechanics. It also allows
to give a formula for $\xs$ for low energies in the periodic
setting only:
\begin{prop}Assume that $N=1$ and that the environment is periodic.

Let $(y,\xi)\in\R^2$ such that $H(y,\xi)<u_{\text{max}}$. Then
$\xs(y,\xi)=0$.

\end{prop}

The organisation of this note is the following : in the second
section, we derive some preliminary results on the long-time
behavior of the system $(Y,\Xi)$ thanks to the ergodic theorem.
Those will be useful in the proof of theorem \ref{thm:general}, to
which is devoted the third section. Eventually, the fourth and
last section is concerned with results in the integrable case,
both in the periodic and the stationary ergodic settings.

\section{Preliminaries}

This section is largely devoted to the study of the long-time
behavior of the Hamiltonian system $(Y,\Xi)$ defined by
\eqref{eq:systhamilt}. First, notice that the Hamiltonian
$H(y,\xi,\om):=\frac 1 2 |\xi|^2 + u(y,\om)$ is constant along the
curves of the system $(Y,\Xi)$, and if $f\in
L^{\infty}(\Omega,\mathcal C^1(\R^N_y\times\R^N_{\xi}))$ is
stationary, then
$$
f\in\mathbb K \ \iff\ f\left(Y(t,y,\xi,\om), \Xi(t,y,\xi,\om),\om
\right)=f(y,\xi,\om)\quad\forall
(y,\xi,\om)\in\R^N\times\R^N\times\Omega.
$$
Indeed, for all $f\in L^{\infty}(\Omega,\mathcal
C^1(\R^N_y\times\R^N_{\xi}))$, we have
$$
\frac{\p }{\p t}f\left(Y(t,y,\xi,\om), \Xi(t,y,\xi,\om),\om
\right)=\left\{H, f  \right\}\left(Y(t,y,\xi,\om),
\Xi(t,y,\xi,\om),\om \right),
$$
where $\left\{H, f  \right\}$ denotes the Poisson bracket of $f$
and $H$, i.e.
$$
\left\{H, f  \right\}(y,\xi,\om)=\xi \cdot\nabla_y f (y,\xi,\om) -
\nabla_y u (y,\xi,\om) \cdot\nabla_{\xi} f (y,\xi,\om).
$$

Let us mention an easily checked property of the trajectories
$(Y,\Xi)$ which will be used extensively in the rest of the
article : for all $(y,z,\xi)\in\R^{3N}$, for all $\om \in\Omega$,
$t\geq 0$,
\begin{gather}
Y(t,y,\xi,\tau_z \om) + z= Y(t,y+z,\xi,\om),\nonumber\\
\Xi(t,y,\xi,\tau_z \om)=\Xi(t,y+z,\xi,\om).
\label{eq:inv_hamiltsyst}\end{gather} In the periodic case, this
invariance entails that the hamiltonian system $(Y,\Xi)$ can be
considered as a dynamical system on the $N$ dimensional torus
$[0,2\pi )^N$. In this periodic setting, it is somewhat natural to
introduce the semi-group of transformations $(\mathcal T_t)_{t\geq
0}$ on $[0,2\pi )^N\times \R^N$ given by
$$
\mathcal T_t(y,\xi)=(Y(t,y,\xi), \Xi(t,y,\xi)),\quad y\in[0,2\pi
)^N,\ \xi\in\R^N.
$$
According to Liouville's theorem, this semi-group preserves the
Lebesgue measure on $[0,2\pi )^N\times \R^N$; moreover, we can
construct a family of finite invariant measures on $[0,2\pi
)^N\times \R^N$ by setting $m_c(y,\xi)=\mathbf 1_{H(y,\xi)\leq
c}\:dy\:d\xi$ for $c>0$ (remember that the Hamiltonian is constant
along the hamiltonian curves). This construction is the root of
the ergodic theorem (see corollary \ref{cor:ergo}), and thus of
the study of the long-time behavior of the system $(Y,\Xi)$.

In the stationary ergodic setting, this construction can be
generalized as follows : we define the transformation
$T_t:\R^N_{\xi}\times \Omega \to\R^N_{\xi}\times \Omega$ by
$$
T_t(\xi, \om)=\left(\Xi(t,0,\xi,\om),\tau_{Y(t,0,\xi,\om)}\om
\right)
$$
together with the family of measures $$\mu_c:=\mathbf 1_{\mathcal
H(\xi,\om)\leq c}\; d\xi\; dP(\om)$$ where $\mathcal
H(\xi,\om):=\frac 1 2 |\xi|^2 + u(0,\om)$. It is obvious that for
all $c\in(0,\infty)$, $\mu_c$ is a finite measure on
$\R^N_{\xi}\times \Omega$.

Notice that the ``good'' generalization to the stationary ergodic
setting of the semi-group $(\mathcal T_t)$ is a semi-group which
acts on $\R^N_{\xi}\times \Omega$ rather than $\R^N_{y}\times
\R^N_{\xi}$. Thanks to the group of transformations
$(\tau_x)_{x\in\R^N}$, the transformations in $\Omega$ can result
in transformations in $\R^N_y$, but the definition chosen here
allows us to define a family of finite invariant measures, whereas
such a construction is rather difficult if one tries to define a
semi-group acting on $\R^N_y\times\R^N_{\xi}$. This will be
fundamental in the rest of the proof.

\begin{lem}
$(T_t)_{t\geq 0}$ is a semi-group on $\R^N_{\xi}\times \Omega$ and
preserves the family of measures $\mu_c$.

\label{lem:Tt}
\end{lem}

\begin{proof}

Let us first prove the semi-group property : let $t,s\in[0,\infty)$,
and $(\xi,\om)\in\R^N\times \Omega$; then
\begin{eqnarray*}
T_t\circ T_s (\xi,\om)&=&T_t
\left(\Xi(s,0,\xi,\om),\tau_{Y(s,0,\xi,\om)}\om \right)\\
&=&\left(\Xi(t,0, \Xi(s,0,\xi,\om),\tau_{Y(s,0,\xi,\om)}\om) ,
\om'\right)
\end{eqnarray*}
and using the properties \eqref{eq:inv_hamiltsyst} we deduce
\begin{eqnarray*}
\Xi(t,0, \Xi(s,0,\xi,\om),\tau_{Y(s,0,\xi,\om)}\om)
&=&\Xi(t,Y(s,0,\xi,\om), \Xi(s,0,\xi,\om),\om),\\
&=&\Xi(t+s, 0,\xi, \om)
\end{eqnarray*}
and
\begin{eqnarray*}
\om'&=&\tau_{Y(t,0,\Xi(s,0,\xi,\om),\tau_{Y(s,0,\xi,\om)}\om)}\tau_{Y(s,0,\xi,\om)}\om\\
&=&\tau_{Y(t,0,\Xi(s,0,\xi,\om),\tau_{Y(s,0,\xi,\om)}\om)
+Y(s,0,\xi,\om) }\om\\
&=&\tau_{Y(t,Y(s,0,\xi,\om),\Xi(s,0,\xi,\om),\om)}\om \\
&=&\tau_{Y(t+s,0,\xi,\om)}\om
\end{eqnarray*}
Thus $$T_t\circ T_s (\xi,\om)=\left( \Xi(t+s, 0,\xi,
\om),\tau_{Y(t+s,0,\xi,\om)}\om \right)=T_{t+s}(\xi,\om).$$

Since it is obvious that $T_0=\Id$, $(T_t)_{t\geq 0}$ is a
semi-group on $\R^N\times\Omega$.

We now have to check the invariance property; let $F\in
L^1(\R^N\times\Omega; \mu_c)$ arbitrary. We set
$f(y,\xi,\om):=F(\xi,\tau_y \om)$ for
$(y,\xi,\om)\in\R^N_y\times\R^N_{\xi}\times \Omega$, and we
compute
$$
\int_{\R^N\times \Omega}F(T_t(\xi,\om))\; d\mu_c(\xi,\om) =
E\left[\int_{R^N}f(Y(t,0,\xi,\om), \Xi(t,0,\xi,\om),\om) \mathbf
1_{H(Y(t,0,\xi,\om),\Xi(t,0,\xi,\om),\om)\leq c}\:d\xi\right]
$$

Since the probability measure $P$ is invariant by the group of
transformation $\tau_y$, and
$$
f(Y(t,y,\xi,\om), \Xi(t,y,\xi,\om),\om)=f(Y(t,0,\xi,\tau_y\om),
\Xi(t,0,\xi,\tau_y\om),\tau_y\om),
$$
we have, for all $y\in\R^N$
\begin{multline*}E\left[f(Y(t,0,\xi,\om),
\Xi(t,0,\xi,\om),\om)1_{H(Y(t,0,\xi,\om),\Xi(t,0,\xi,\om),\om)\leq
c}\right]=\\=E\left[f(Y(t,y,\xi,\om),
\Xi(t,y,\xi,\om),\om)1_{H(Y(t,y,\xi,\om),\Xi(t,y,\xi,\om),\om)\leq
c}\right].\end{multline*}

Take an arbitrary function $\phi\in L^1(\R^N_y)$, and write
\begin{eqnarray*}
&&\int_{\R^N\times \Omega}F(T_t(\xi,\om))\; d\mu_c(\xi,\om)\\
&=&E\left[\int_{\R^{2N}}dy\:d\xi \phi(y)f(Y(t,y,\xi,\om),
\Xi(t,y,\xi,\om),\om)1_{H(Y(t,y,\xi,\om),\Xi(t,y,\xi,\om),\om)\leq
c}\right]
\end{eqnarray*}
We change variables in the integral in $(y,\xi)$ by setting
$(x,v)=(Y(t,y,\xi,\om), \Xi(t,y,\xi,\om))$; according to Liouville's
theorem, the jacobian of this change of variables is equal to 1, and
$$
(x,v)=(Y(t,y,\xi,\om), \Xi(t,y,\xi,\om))\iff (y,\xi)=(X(t,x,v,\om),
V(t,x,v,\om)),
$$
where $(X,V)$ is a solution of the Hamiltonian system
$$
\left\{\begin{array}{l} \dot{X}=V,\\\dot{V}=-\nabla u
(X,\om),\\(X,V)(t=0,x,v)=(x,v).\end{array}\right.
$$

Observe that in the present case, we have simply
$$
X(t,x,v,\om)=Y(t,x,-v,\om),
$$
so that $(X,V)$ satisfies relations \eqref{eq:inv_hamiltsyst}.

Hence\begin{eqnarray*}&& \int_{\R^{2N}}dy\:d\xi
\phi(y)f(Y(t,y,\xi,\om),
\Xi(t,y,\xi,\om),\om)1_{H(Y(t,y,\xi,\om),\Xi(t,y,\xi,\om),\om)\leq
c}\\
&=&\int_{\R^{2N}}dx\:dv \phi(X(t,x,v,\om))f(x,
v,\om)1_{H(x,v,\om)\leq c}\\
&=&\int_{\R^{2N}}dx\:dv \phi(X(t,0,v,\tau_x\om) + x)F(
v,\tau_x\om)1_{\mathcal H(v,\tau_x\om)\leq c}
\end{eqnarray*}
so that
\begin{eqnarray*}
&&\int_{\R^N\times \Omega}F(T_t(\xi,\om))\; d\mu_c(\xi,\om)\\
&=& E\left[\int_{\R^{2N}}dx\:dv \phi(X(t,0,v,\tau_x\om) + x)F(
v,\tau_x\om)1_{\mathcal H(v,\tau_x\om)\leq c} \right]\\
&=&E\left[\int_{\R^{2N}}dx\:dv \phi(X(t,0,v,\om) + x)F(
v,\om)1_{\mathcal H(v,\om)\leq c} \right]\\
&=&E\left[\int_{\R^{N}}\:dv \left(\int_{\R^N}\phi(X(t,0,v,\om) +
x)\;
dx\right)F( v,\om)1_{\mathcal H(v,\om)\leq c} \right]\\
&=&E\left[\int_{\R^{N}}\:dv F( v,\om)1_{\mathcal H(v,\om)\leq c}
\right]=\int_{\R^N\times\Omega}F \:d \mu_c
\end{eqnarray*}
since the integral of $\phi$ is equal to 1.

\end{proof}
The following corollary is an immediate consequence of Birkhoff's
ergodic theorem:

\begin{corol}
Let $F\in L^1(\R^N\times\Omega; \mu_c)$. There exists a function
$\bar{F}\in L^1(\R^N\times\Omega; \mu_c)$ such that as $T\to
\infty$,
$$
\frac 1 T \int_0^T F(T_t(\xi,\om))\;dt \to \bar{F} (\xi,\om)
$$
a.e. on $\R^N\times\Omega$ and in $L^1(\mu_c)$. Moreover,
$\bar{F}$ is invariant by $T_t$ for all $t>0$, and \be
\int_{\R^N\times \Omega}F\: d\mu_c = \int_{\R^N\times
\Omega}\bar{F}\: d\mu_c. \label{eq:invFbarF}\ee

Additionally, if $\bar{f}=\bar{f}(y,\xi,\om)$ is the stationary
function associated to $\bar{F}$, that is,
$\bar{f}(y,\xi,\om)=\bar{F}(\xi,\tau_y\om)$, then $\bar{f}$ is
invariant by the hamiltonian flow $(Y,\Xi)$; precisely, for all
$(y,\xi,\om)\in\R^{2N}\times\Omega$, $t>0$
$$
\bar{f}(Y(t,y,\xi,\om),\Xi(t,y,\xi,\om), \om)=\bar{f}(y,\xi,\om).
$$
\label{cor:ergo}
\end{corol}
\begin{proof}
We only have to prove the invariance of $\bar{f}$ by the Hamiltonian
flow; first, for $y=0$, we have
\begin{eqnarray*}
\bar{f}(Y(t,0,\xi,\om),\Xi(t,0,\xi,\om),
\om)&=&\bar{F}(\Xi(t,0,\xi,\om),
\tau_{Y(t,0,\xi,\om)}\om)=\bar{F}(T_t(\xi,\om))\\
&=&\bar{F}(\xi,\om)=\bar{f}(0,\xi,\om)
\end{eqnarray*}
and the property is proved when $y=0$.

For $y\in\R^N$ arbitrary,
\begin{eqnarray*}
\bar{f}(Y(t,y,\xi,\om),\Xi(t,y,\xi,\om),
\om)&=&\bar{f}(Y(t,0,\xi,\tau_y\om) + y,\Xi(t,0,\xi,\tau_y\om),
\om)\\
&=&\bar{f}(Y(t,0,\xi,\tau_y\om) ,\Xi(t,0,\xi,\tau_y\om),
\tau_y\om)\\
&=&\bar{f}(0,\xi,\tau_y\om)=\bar{f}(y,\xi,\om)
\end{eqnarray*}
according to the result in the case $y=0$.

\end{proof}

\begin{rem}
We mention here an important but easy consequence of the relations
\eqref{eq:inv_hamiltsyst} and the invariance of the measure $P$
w.r.t. $\tau_y$, $y\in\R^N$ : for any stationary function
$f=f(y,\xi,\om)=F(\xi,\tau_y\om)$, $F\in
L^{\infty}(\R^N\times\Omega)$, we have
$$
E[f(Y(t,y,\xi,\cdot),
\Xi(t,y,\xi,\cdot),\cdot)]=E[F(T_t(\xi,\cdot))]
$$
for all $t>0$, $y,\xi\in\R^N$; in particular, the left-hand side of
the above equality does not depend on $y$.

This property was used in the proof of lemma \ref{lem:Tt}
\label{rem:statio}
\end{rem}

\begin{rem}
Let us precise a little what happens when the function $F\in
L^1_{\text{loc}}(\R^N_{\xi},L^1(\Omega))$. In that case, $F\in
L^1(\R^N_{\xi}\times\Omega; \mu_c)$ for all $c>0$. Consequently,
for any $c>0$, we can define the function $\bar{F}_c$ associated
to $F$ by corollary \ref{cor:ergo}.

It is then easily proved that for any $0<c<c'$,
$\bar{F}_c=\bar{F}_{c'}$, $\mu_c$-almost everywhere. Setting
$A_n=\{(\xi,\om)\in\supp \mu_n;\ \bar{F}_n(y,\xi)\neq
\bar{F}_{n+1}(y,\xi)\} $, and $A=\cup_{n=0}^{\infty} A_n$, we see
that $\mu_c(A)=0$ for all $c>0$. Moreover, for all
$(\xi,\om)\in\R^N\times\Omega\setminus A$, for all integers $k,l$
such that $(\xi,\om)\in \supp \mu_k\cap \supp \mu_l$, we have
$\bar{F}_k(\xi,\om)=\bar{F}_l(\xi,\om)$. We can thus define a
function $\bar{F}(\xi,\om)$ on $\R^N\times\Omega\setminus A$ by
$$
\bar{F}(\xi,\om)=\bar{F}_n(\xi,\om)\quad \text{for any $n\in\N$
such that $(\xi,\om)\in \supp \mu_n$}
$$
We then now that \be \frac 1 T \int_0^T F(T_t(\xi,\om))\:dt \to
\bar{F}(\xi,\om) \label{conv_ergo}\ee as $T\to \infty$, and the
convergence holds in $L^1(\mu_c)$ for all $c>0$, and $\mu_n$
almost everywhere for $n\in\N$. Eventually, setting
$$
B:=\{(\xi,\om)\in\R^N\times\Omega\setminus A ;\ \frac 1 T \int_0^T
F(T_t(\xi,\om))\:dt \ \text{does not converge towards}\
\bar{F}(\xi,\om)\}
$$
it is easily proved that $\mu_c(B)=0$ for all $c>0$ (the equality
is true for $c\in\N$, and is then deduced for $c>0$ arbitrary
because the family of measures $(\mu_c)$ is increasing in $c$).

Eventually, we have found a function $\bar{F}\in
L^1_{\text{loc}}(\R^N, L^1(\Omega))$, independent of $c$, such
that \eqref{conv_ergo} holds in $L^1(\mu_c)$ and $\mu_c$-almost
everywhere for all $c>0$.

\label{rem:proj_indt_c}
\end{rem}

\begin{rem}
The construction above allows us to make more precise what we mean
by projection $P$: let $f=f(y,\xi,\om)$ be a stationary function,
$f\in L^{\infty}(\R^N_y, L^1_{\text{loc}}(\R^N_{\xi},
L^1(\Omega)))$, and set $F(\xi,\om)=f(0,\xi,\om)\in
L^1_{\text{loc}}(\R^N_{\xi},L^1(\Omega))$. We can then associate
to $F$ a function $\bar{F}\in
L^1_{\text{loc}}(\R^N_{\xi},L^1(\Omega))$ such that
\eqref{conv_ergo} holds in $L^1(\mu_c)$ for all $c$ (see remark
\ref{rem:proj_indt_c}). We set
$$
P(f)(y,\xi,\om):=\bar{F}(\xi,\tau_y\om).
$$
It follows from corollary \ref{cor:ergo} that $P(f)$ is invariant
by the hamiltonian flow \eqref{eq:systhamilt}, and thus satisfies
the constraint equation. From now on, we take this definition for
the projection $P$, instead of the one given in the introduction.
Notice that, for all $y\in\R^N$ and $\mu_c$-almost everywhere,
\begin{eqnarray*}
P(f)(y,\xi,\om)&=&\lim_{T\to\infty} \frac 1 T \int_0^T
F(T_t(\xi,\tau_y\om))\:dt\\
&=&\lim_{T\to\infty} \frac 1 T \int_0^T
f\left(Y(t,0,\xi,\tau_y\om),\Xi(t,0,\xi,\tau_y\om),\tau_y\om\right)\:dt\\
&=&\lim_{T\to\infty} \frac 1 T \int_0^T
f\left(Y(t,y,\xi,\om),\Xi(t,y,\xi,\om),\om\right)\:dt
\end{eqnarray*}

And we also give a more precise definition of $\xs(y,\xi,\om)$ :
let$$
\hat{\xi}:\begin{array}{ccc} \R^N\times \Omega &\to&\R^N\\
                                (\xi,\om)&\mapsto&\xi
\end{array}.
$$
Then
$$
\xs(y,\xi,\om)=P(\hat{\xi})(y,\xi,\om)=\lim_{T\to\infty} \frac 1 T
\int_0^T \Xi(t,y,\xi,\om)\:dt
$$
a.e. and in $L^1(\mu_c)$ for all $0<c<\infty$.

\label{rem:P_equiv_ergo}
\end{rem}

Eventually, we mention here a property that will be used in the
proof of the theorem; with the same notations as above, let
$$
\phi(\tau,y,\xi,\om)=F\left( T_{\tau}(\xi,\tau_y\om)\right).
$$
Then $\phi$ is a solution of the evolution equation
$$
\p_{\tau}\phi + \xi\cdot \nabla_y \phi -  \nabla_y u \cdot
\nabla_{\xi} \phi=0,
$$
with initial data $\phi(\tau=0,y,\xi,\om)=f(y,\xi,\om)=F\left(
\xi,\tau_y\om\right)$.
\section{The general $N$-dimensional case}

This section is devoted to the proof of theorem \ref{thm:general}.
The proof is divided in three steps : first, we study the case of
an initial data which does not depend on $x$, then the case when
the initial data only depends on $x$ (and not on $y,\xi,\om$), and
eventually, we treat the general case.

\subsection{First case : $f_0$ does not depend on $x$}

Here, we assume that $f_0=f_0(y,\xi,\om)\in
L^1_{\text{loc}}(\R^N_{\xi}; L^{\infty}(\R^N_{y}\times \Omega)$.
Recall that $f_0$ is stationary, i.e.
$f_0(y+z,\xi,\om)=f_0(y,\xi,\tau_z\om)$ a.s. in $\om$, for all
$(y,z,\xi)\in\R^{3N}$. In the rest of the subsection, we set
$$
F_0(\xi,\om):=f_0(0,\xi,\om)
$$
and
$$
\bar{F}_0(\xi,\om):=\lim_{T\to \infty} \frac{1}{T}\int_0^T
F(T_t(\xi,\om))\:dt,\quad
\bar{f}_0(y,\xi,\om)=\bar{F}_0(\xi,\tau_y\om).
$$
Notice that $F\in L^1_{\text{loc}}(\R^N_{\xi};
L^{\infty}(\Omega))$, and thus $F\in L^1(\R^N\times\Omega; \mu_c)$
for all $c>0$.

In that case, \begin{eqnarray*}
\fe(t,x,\xi,\om)&=&f_0\left(Y\left(\frac{t}{\e},\frac{x}{\e},\xi,\om
\right),\Xi\left(\frac{t}{\e},\frac{x}{\e},\xi,\om
\right),\om\right)\\
&=&f_0\left(Y\left(\frac{t}{\e},0,\xi,\tau_{\frac{x}{\e}}\om
\right),\Xi\left(\frac{t}{\e},0,\xi,\tau_{\frac{x}{\e}}\om
\right),\tau_{\frac{x}{\e}}\om\right)\\
&=&F_0\left(T_{\frac{t}{\e}}\left(\xi,
\tau_{\frac{x}{\e}}\om\right)\right)\\
&=&\bar{f}_0\left( \frac{x}{\e},\xi,\om\right) +
\left\{F_0\left(T_{\frac{t}{\e}}\left(\xi,
\tau_{\frac{x}{\e}}\om\right)\right)-
\bar{F}_0\left(\xi,\tau_{\frac{x}{\e}}\om\right)\right\}
\end{eqnarray*}

In accordance with theorem \ref{thm:general}, we set
$$
g(\tau,y,\xi,\om)=\left(F_0
-\bar{F}_0\right)\left(T_{\tau}\left(\xi,
\tau_{y}\om\right)\right),
$$
and $\re=0$. Then $g$ satisfies the microscopic evolution equation
\eqref{eq:evol_micro} thanks to the remark at the end of the
preceding section. Moreover, $g(\tau)\in\K^{\bot}$ by definition
of $\K^{\bot}$ and because  $P\left(F_0\left(T_{\tau}\left(\xi,
\tau_{y}\om\right)\right)\right)=\bar{F}_0(\xi,\tau_y\om)$.

It only remains to check that
$$
\int_0^T g\left(\frac{t}{\e}, \frac{x}{\e},\xi,\om\right)\:dt\to
0\quad \text{as}\ \e\to 0
$$
in $L^1_{\text{loc}}(\R^N_x,L^1(\R^N\times\Omega, \mu_c))$ for all
$T>0$ and $c>0$.

The invariance of the measure $P$ with respect to the group of
transformations $(\tau_x)_{x\in\R^N}$ (see remark
\ref{rem:statio}) entails that
\begin{eqnarray*}
&&\int_{\Omega\times\R^N_{\xi}}
\left|\frac{1}{\;\frac{T}{\e}\;}\int_0^{\frac{T}{\e}}f_0\left(Y\left(t,\frac{x}{\e},\xi,\om
\right),\Xi\left(t,\frac{x}{\e},\xi,\om \right),\om\right)\:dt
-\bar{f}_0\left(\frac{x}{\e},\xi,\om \right)\right|\:
\:d\mu_c(\xi,\om)\\
&=&
\int_{\Omega\times\R^N_{\xi}}\left|\frac{1}{\;\frac{T}{\e}\;}\int_0^{\frac{T}{\e}}F_0\left(T_t(\xi,\om)
\right)\:dt -\bar{F}_0\left(\xi,\om \right)\right|
\:d\mu_c(\xi,\om)
\end{eqnarray*}
and the term above goes to 0 as $\e\to 0$ according to corollary
\ref{cor:ergo} and is independent of $x\in\R^N$. There remains to
check that $\bar{f}_0=P(f_0)$. This follows directly from remark
\ref{rem:P_equiv_ergo}. Thus theorem \ref{thm:general} is proved
in the case when $f_0$ does not depend on the macroscopic variable
$x$. \vskip1mm

The following remark will prove to be useful when treating the
general case :
\begin{rem}
If $f_0\in L^{\infty}$, then for any function $a\in
L^{\infty}((0,\infty)\times\R^N_y\times\R^N_{\xi}\times\Omega)$,
stationary in $y$, we have
$$
\int_0^T a\left(t,\frac{x}{\e},\xi,\om\right)g\left(\frac{t}{\e},
\frac{x}{\e},\xi,\om\right)\:dt\to 0\quad \text{as}\ \e\to 0
$$
in $L^1_{\text{loc}}(\R^N_x,L^1(\R^N\times\Omega, \mu_c))$ for all
$T>0$ and $c>0$.

Indeed, prove the property first for $a=a_1(t)a_2(y,\xi,\om)$,
with $a_1$, $a_2\in L^{\infty}$. For $a$ arbitrary, take a
sequence $a_{\delta}$ with $\delta>0$, converging to $a$ in
$L^1_{\text{loc}}$, and such that
$$
a_{\delta}=\sum_{k=0}^{n_{\delta}}a_1^{\delta}(t)a_2^{\delta}(y,\xi,\om).
$$
with $a_1^{\delta},a_2^{\delta}$ in $L^{\infty}$. The property is
known for $a_{\delta}$, and it is thus easily deduced for $a$.

\label{rem:wkcv_0}
\end{rem}

\subsection{Second case : $f_0=f_0(x)$ }

Unlike the preceding subsection, we now focus on the case when
$f_0$ only depends on the macroscopic variable $x$. In order to
simplify the analysis, we assume that $f_0 \in
 W^{1,\infty}(\R^N_x)$ (the case
when $f_0$ is not smooth in $x$ will be treated in the next
subsection). In that case,
$$
\fe(t,x,\xi,\om)=f_0\left(\e Y\left(
\frac{t}{\e},\frac{x}{\e},\xi,\om\right)\right).
$$
Hence we have to investigate the behavior as $\e\to 0$ of
$$
\e Y\left( \frac{t}{\e},\frac{x}{\e},\xi,\om\right).
$$
We prove the following
\begin{lem}
Let $T>0$ arbitrary. As $\e $ vanishes,
$$
\e Y\left( \frac{t}{\e},\frac{x}{\e},\xi,\om\right)- x +
t\xi^{\sharp}\left(\frac{x}{\e},\xi,\om\right) \to 0
$$
in $L^{\infty}((0,T)\times\R^N_x;
L^1(\R^N_{\xi}\times\Omega,\mu_c))$.

\end{lem}
\begin{proof}
Let us write, for $t>0$
\begin{eqnarray*}
\e Y\left( \frac{t}{\e},\frac{x}{\e},\xi,\om\right)- x +t
\xi^{\sharp}\left(\frac{x}{\e},\xi,\om\right)&=&\e
\int_0^{\frac{t}{\e}}\dot{Y}\left(
s,\frac{x}{\e},\xi,\om\right)\:ds +
t\xi^{\sharp}\left(\frac{x}{\e},\xi,\om\right)\\
&=&-t \frac{\e}{t}\int_0^{\frac{t}{\e}}\Xi\left(
s,\frac{x}{\e},\xi,\om\right)\:ds +
t\xi^{\sharp}\left(\frac{x}{\e},\xi,\om\right)\\
&=&-t\left\{\frac{\e}{t}\int_0^{\frac{t}{\e}}\hat{\xi}\left(T_s(\xi,\tau_{\frac{x}{\e}}\om)
\right)\:ds -\xs\left(\frac{x}{\e},\xi,\om\right)\right\}
\end{eqnarray*}
Let $0<\alpha<T$ arbitrary. For $\alpha\leq t\leq T$, we have
\begin{eqnarray*}
&&\int_{\R^N_{\xi}\times\Omega}\left|\e Y\left(
\frac{t}{\e},\frac{x}{\e},\xi,\om\right)- x +t
\xi^{\sharp}\left(\frac{x}{\e},\xi,\om\right)\right|\:d\mu_c
(\xi,\om)\\&=&t
\int_{\R^N_{\xi}\times\Omega}\left|\frac{\e}{t}\int_0^{\frac{t}{\e}}\hat{\xi}\left(T_s(\xi,\om)
\right)\:ds -\xs\left(0,\xi,\om\right)\right|\:d\mu_c (\xi,\om)\\
&\leq & T \sup_{\tau\geq \frac{\alpha}{\e}}\left|\left|
\frac{1}{\tau}\int_0^{\tau}\hat{\xi}\left(T_s(\xi,\om) \right)\:ds
-\xs\left(0,\xi,\om\right)\right|\right|_{L^1(\R^N\times\Omega,\mu_c)}
\end{eqnarray*}
and the upper-bound vanishes as $\e\to 0$ for any $\alpha>0$
thanks to corollary \ref{cor:ergo}. Notice that the upper-bound
does not depend on $x$, hence the convergence holds in
$L^{\infty}(\R^N_x; L^1(\mu_c))$.

We now have to investigate what happens when $t$ is close to 0;
notice that
$$
\sup_{x\in\R^N} ||\xi^{\sharp}\left(\frac{x}{\e},\xi,\om\right)
||_{L^1(\R^N\times\Omega, \mu_c)}\leq C_0
$$
where the constant $C_0$ only depends on $N$ and $c$. Similarly,
for all $t\geq 0$,
$$
\sup_{x\in\R^N}
\left|\left|\hat{\xi}\left(T_s(\xi,\tau_{\frac{x}{\e}})\right)\right|\right|_{L^1(\R^N\times\Omega,
\mu_c)}\leq C_0.
$$

Hence, if $0\leq t\leq \alpha$, we have
$$
\sup_{x\in \R^N}\int_{\R^N_{\xi}\times\Omega}\left|\e Y\left(
\frac{t}{\e},\frac{x}{\e},\xi,\om\right)- x +
\xi^{\sharp}\left(\frac{x}{\e},\xi,\om\right)\right|\:d\mu_c
(\xi,\om)\leq 2\alpha C_0.
$$
Eventually,
\begin{multline*}
\left| \left|\e Y\left( \frac{t}{\e},\frac{x}{\e},\xi,\om\right)-
x + t\xi^{\sharp}\left(\frac{x}{\e},\xi,\om\right)\right|
\right|_{L^{\infty}((0,T)\times\R^N; L^1(\mu_c))} \leq\\\leq
\inf_{0<\alpha<T}\left\{2C_0 \alpha +T \sup_{\tau\geq
\frac{\alpha}{\e}}\left|\left|
\frac{1}{\tau}\int_0^{\tau}\hat{\xi}\left(T_s(\xi,\om )\right)\:ds
-\xs\left(0,\xi,\om\right)\right|\right|_{L^1(\mu_c)} \right\}
\end{multline*}
and the lemma is proved.
\end{proof}

We easily deduce that theorem \ref{thm:general} is true when
$f_0\in W^{1,\infty}(\R^N)$ with
\begin{gather*}
f(t,x,y,\xi,\om):=f_0(x-t\xs(y,\xi,\om)),\quad g=0,\\
r^{\e}(t,x,\xi,\om):=\fe(t,x,\xi,\om)-f\left(t,x,\frac{x}{\e},\xi,\om\right)
\end{gather*}
and it is easily checked that $f$ satisfies $P(f)=f$,
$f(t=0)=P(f_0)=f_0$ (since $f_0$ is independent of $y$ and $\xi$),
and that $f$ is a solution of the evolution equation
\eqref{eq:evol_homog}.

\subsection{Third case : $f_0$ arbitrary}

We now tackle the case of an arbitrary stationary function $f_0\in
L^1_{\text{loc}}(\R^N_x\times\R^N_{\xi},
L^{\infty}(\R^N_y\times\Omega))$. We begin with the case when
$$
f_0(x,y,\xi,\om)=a(x) b(y,\xi,\om),
$$
with $a\in  W^{1,\infty}(\R^N)$ and $b\in
L^1_{\text{loc}}(\R^N_{\xi},L^{\infty}(\R^N_y\times\Omega))\cap
L^{\infty}(\R^N_y\times\R^N_{\xi}\times\Omega)$, $b$ stationary.
This case follows directly from the two first subsections. Indeed,
let
$$
f(t,x,y,\xi,\om)=a(x-t\xs(y,\xi,\om))\:P(b)(y,\xi,\om),
$$
and
$$
g(t,x; \tau,y,\xi,\om)=a(x-t\xs(y,\xi,\om))\:\left(b-P(b)
\right)(T_{\tau}(y,\xi,\om)).
$$
It is already known that $f$ and $g$ satisfy
\eqref{eq:evol_homog}, that $f(t,x,\cdot)\in \mathbb K$, and that
$g$ satisfies \eqref{eq:evol_micro} thanks to the preceding
subsections and the fact that $\xs(y,\xi,\om)$ is invariant by the
Hamiltonian flow $(Y,\Xi)$. Notice that it is capital here that
the coefficient $\xs(y,\xi,\om)$ in the transport equation
\eqref{eq:evol_homog} belongs to $\K$.

There remains to check that $g(t,x;\tau,\cdot)\in \K^{\bot}$, that
the remainder $\re$ goes to 0 strongly in $L^1_{\text{loc}}$ and
that $g(t,x; t/\e,x/\e,\xi,\om)$ goes weakly to 0 in the sense of
theorem \ref{thm:general}. First, notice that
$a(x-t\xs(y,\xi,\om))\in \K $ and $\left(b-P(b)
\right)(T_{\tau}(y,\xi,\om)\in \K^{\bot}$. Thus, $
a(x-t\xs)P(b)=P(a(x-t\xs)b)$ almost everywhere (because
$a(x-t\xs(0,\xi,\om))$ is invariant by the semi-group $T_{\tau}$),
and consequently
$$
g(t,x; \tau,y,\xi,\om)=\left[a(x-t\xs) b - P\left(a(x-t\xs) b
\right) \right](T_{\tau}(\xi,\tau_y\om)).
$$
Hence $g(t,x;\tau)\in \K^{\bot}$ a.e.

Then, setting
$$
\re(t,x,\xi,\om)=\fe(t,x,\xi,\om) -
f\left(t,x,\frac{x}{\e},\xi,\om \right) - g\left(t,x;
\frac{t}{\e},\frac{x}{\e},\xi,\om \right),
$$
we have to prove that $\re$ goes to 0 strongly in
$L^1_{\text{loc}}$. We compute the difference
\begin{eqnarray*}
&&\fe(t,x,\xi,\om ) - f\left(t,x,\frac{x}{\e},\xi,\om\right) -
g\left(t,x; \frac{t}{\e},\frac{x}{\e},\xi,\om \right)
\\&=&a\left(\e Y\left(
\frac{t}{\e},\frac{x}{\e},\xi,\om\right)\right) b\left(Y\left(
\frac{t}{\e},\frac{x}{\e},\xi,\om\right),\Xi\left(
\frac{t}{\e},\frac{x}{\e},\xi,\om\right),\om \right) \\
&&- a\left(x-t\xs\left(\frac{x}{\e},\xi,\om\right)
\right)\:P(b)\left(\frac{x}{\e},\xi,\om\right)\\
&&-a\left(x-t\xs\left(\frac{x}{\e},\xi,\om\right)\right) \left[ b-
P(b)\right]\left(Y\left(
\frac{t}{\e},\frac{x}{\e},\xi,\om\right),\Xi\left(
\frac{t}{\e},\frac{x}{\e},\xi,\om\right),\om \right)\\
&=&\left[a\left(\e Y\left(
\frac{t}{\e},\frac{x}{\e},\xi,\om\right)\right)-
a\left(x-t\xs\left(\frac{x}{\e},\xi,\om\right)
\right)\right]\:b\left(Y\left(
\frac{t}{\e},\frac{x}{\e},\xi,\om\right),\Xi\left(
\frac{t}{\e},\frac{x}{\e},\xi,\om\right),\om \right)
\end{eqnarray*}

The right-hand side of the above equality is bounded by
$$
||a||_{W^{1,\infty}} ||b||_{L^{\infty}}\left|\e Y\left(
\frac{t}{\e},\frac{x}{\e},\xi,\om\right)-x + t
\xs\left(\frac{x}{\e},\xi,\om\right)\right|
$$
and thus converges to 0 as $\e\to 0$ in
$L^{\infty}((0,T)\times\R^N_x;
L^1(\R^N_{\xi}\times\Omega,\mu_c)))$ according to the second
subsection.

Moreover, it is easily proved that as $\e\to 0$,
$$
\int_0^T g\left(t,x; \frac{t}{\e}, \frac{x}{\e},\xi,\om
\right)\;dt\to 0
$$
strongly in $L^1_{\text{loc}}(\R^N_x\times\R^N_{\xi},
L^1(\Omega))$ thanks to remark \ref{rem:wkcv_0}. Hence theorem
\ref{thm:general} is proved in that case.

Now, let $f_0\in L^1_{\text{loc}}(\R^N_x\times\R^N_{\xi},
L^{\infty}(\R^N_y\times\Omega))$ arbitrary, and set
$F_0(x,\xi,\om):=f_0(x,0,\xi,\om)$. Take a sequence of functions
$F_n\in L^1(\R^N_x\times\R^N_{\xi}\times\Omega)$ such that
\begin{itemize}
\item $F_n\to F_0$ as $n \to \infty$ in
$L^1_{\text{loc}}(\R^N_x\times\R^N_{\xi}\times\Omega)$;

\item For all $n\in\N$, there exist functions $a_k^n\in L^1\cap
W^{1,\infty}(\R^N)$, $b_k^n\in L^1\cap
L^{\infty}(\R^N_{\xi}\times\Omega)$, $1\leq k\leq n$ such that
$$
F_n(x,\xi,\om)=\sum_{k=1}^n a_k^n(x)\: b_k^n(\xi,\om)\quad
\text{a.e}.
$$
\end{itemize}
Let $\fe_n$ be the solution of \eqref{eq:kinhomog} with initial
data $F_n\left(x,\xi,\tau_{\frac{x}{\e}} \om\right)$, and let
$f_n=f_n(t,x,y,\xi,\om)$, $g_n=g_n(t,x;\tau,y,\xi,\om)$ be the
functions associated to $\fe_n$ by theorem \ref{thm:general} for
all $n$.

Let $f(t,x,y,\xi,\om)$, $g(t,x;\tau,y,\xi,\om)$ be the solutions
of the system
\begin{gather*}
P(f)=f,\quad P(g)=0,\\
\p_t \left(\begin{array}{c}f\\g\end{array}\right) + \xs(y,\xi,\om)\cdot \nabla_x \left(\begin{array}{c}f\\g\end{array}\right)=0,\\
\p_{\tau} g + \xi\cdot\nabla_y g - \nabla_y
u(y,\om)\cdot\nabla_{\xi} g=0,\\ f(t=0)=P(f_0),\quad
g(t=0,x;\tau=0,y,\xi,\om)=\left[f_0 -P(f_0)\right](x,y,\xi,\om).
\end{gather*}

We have already proved that $f_n$, $g_n$ satisfy the above system.
We denote by $\bar{F}_0$, $\bar{F}_n$, the functions associated to
$F_0$, $F_n$ respectively by corollary \ref{cor:ergo}, so that
$P(f_0)(x,y,\xi,\om)=\bar{F}_0(x,\xi,\tau_y\om)$, and
$f_n(t=0,x,y,\xi,\om)=\bar{F}_n(x,\xi,\tau_y\om)$,
$g_n(t=0,x,\tau=0,y,\xi,\om)=(F_n -\bar{F}_n)(x,\xi,\tau_y\om)$.

We use the following lemma, of which we postpone the proof :
\begin{lem}
Let $\geps$ be a solution of \eqref{eq:kinhomog} with initial data
$g_0\left(x,\frac{x}{\e}, \xi,\om\right)$, and $g_0\in
L^1_{\text{loc}}(\R^N_x\times\R^N_{\xi},
L^{\infty}(\R^N_y\times\Omega))$ stationary. Then for all
$R,R',T>0$, for all $t\in[0,T]$,
$$
\int_{x\in B_R,\ \xi\in
B_{R'}}\left|\geps(t,x,\xi,\om)\right|\:dx\:d\xi\leq
||g_0||_{L^1(K_{T,R,R'}, L^{\infty}(\R^N_y\times\Omega))}
$$
where
$$
K_{T,R,R'}=\left\{(x,\xi)\in\R^N\times\R^N, |x|\leq R + T
\sqrt{{R'}^2 + 2 u_{\text{max}}},\ |\xi|\leq \sqrt{{R'}^2 + 2
u_{\text{max}}} \:\right\}.
$$

Similarly, if $g$ is a solution of \eqref{eq:evol_homog} with
initial data $g_0\in L^1_{\text{loc}}(\R^N_x\times\R^N_{\xi},
L^{\infty}(\R^N_y\times\Omega))$, then
$$
\int_{x\leq R}\left|g(t,x,y,\xi,\om)\right|\:dx\leq \int_{x\leq R
+ T \sqrt{{\xi}^2 + 2
u_{\text{max}}}}\left|g_0(x,y,\xi,\om)\right|\:dx
$$

\label{lem:finite_prop_speed}

\end{lem}

Consequently, with $C_{R,T,\xi}:=\{x\in\R^N, \;|x|\leq R + T
\sqrt{{\xi}^2 + 2 u_{\text{max}}}\}$, we have
\begin{eqnarray*}
\left|\left|f(t,\cdot,y,\xi,\om)-f_n(t,\cdot,y,\xi,\om)\right|\right|_{L^1(B_R)}&\leq
&\left|\left|\bar{F}_0(\cdot,\xi,\tau_y\om)-\bar{F}_n(\cdot,\xi,\tau_y\om)\right|\right|_{L^1(C_{R,T,\xi})}\\
&\leq &\left|\left|\overline{|F_0-F_n|}(\cdot,\xi,\tau_y\om)\right|\right|_{L^1(C_{R,T,\xi})}\\
\left|\left|f(t,x,y,\xi,\om)-f_n(t,x,y,\xi,\om)\right|\right|_{L^1(B_R\times\R^N_{\xi}\times\Omega,
dx\:d\mu_c(\xi,\om))}&\leq&\left|\left|\overline{|F_0-F_n|}(x,\xi,\om)\right|\right|_{L^1(C_{R,T,\sqrt{2c}}\times\R^N_{\xi}\times\Omega,
dxd\mu_c(\xi,\om))}\\
&\leq
&\left|\left|F_0-F_n(x,\xi,\om)\right|\right|_{L^1(C_{R,T,\sqrt{2c}}\times\R^N_{\xi}\times\Omega,
dxd\mu_c(\xi,\om))}.
\end{eqnarray*}
In the last inequality, we have used property \eqref{eq:invFbarF}.

And similarly,
\begin{gather*}
\left|\left|\fe(t,x,\xi,\om)-\fe_n(t,x,\xi,\om)\right|\right|_{L^1(B_R\times
B_{R'}\times\Omega)}\leq
\left|\left|F_0-F_n\right|\right|_{L^1(K_{T,R,R'},
L^{\infty}(\Omega))},\\
\left|\left|g(t,x;\tau,y,\xi,\om)-g_n(t,x;\tau,y,\xi,\om)\right|\right|_{L^1(B_R\times\R^N_{\xi}\times\Omega,
dx\:d\mu_c(\xi,\om))} \leq
2\left|\left|F_0-F_n(x,\xi,\om)\right|\right|_{L^1(C_{R,T,\sqrt{2c}}\times\R^N_{\xi}\times\Omega,dxd\mu_c(\xi,\om))}.
\end{gather*}

The above inequalities are true for all $t\in[0,T]$ and for all
$\tau \geq 0$.

Set
$$
\re(t,x,\xi,\om):=\fe(t,x,\xi,\om) -
f\left(t,x,\frac{x}{\e},\xi,\om \right)-
g\left(t,x;\frac{t}{\e},\frac{x}{\e},\xi,\om \right).
$$
Then for all $t\in[0,T]$, for all $n\in\N$, setting $c=\frac 1 2
{R'}^2 + u_{\text{max}}$,
\begin{eqnarray*}
||\re(t)||_{L^1(B_R\times B_{R'}\times\Omega)}&\leq
&\left|\left|\fe(t)-\fe_n(t)\right|\right|_{L^1(B_R\times
B_{R'}\times\Omega)}\\
&&+
||f(t)-f_n(t)||_{L^{\infty}(\R^N_y;L^1({B_R}_x\times\R^N_{\xi}\times\Omega,
dx\:d\mu_c(\xi,\om)))}\\
&& +||g(t)-g_n(t)||_{L^{\infty}((0,\infty)_{\tau}\times\R^N_y;
){L^1(B_R\times\R^N_{\xi}\times\Omega, dx\:d\mu_c(\xi,\om))}}\\
&& + ||\re_n(t)||_{L^1(B_R\times B_{R'}\times\Omega)}\\
&\leq &4
\left|\left|F_0-F_n\right|\right|_{L^1(C_{R,T,\sqrt{2c}}\times\R^N_{\xi}\times\Omega,
dxd\mu_c(\xi,\om))} +||\re_n(t)||_{L^1(B_R\times
B_{R'}\times\Omega)}
\end{eqnarray*}
Thus $\re\to 0$ as $\e\to 0$ in $L^{\infty}([0,\infty);
L^1_{\text{loc}}(\R^N_x\times\R^N_{\xi}; L^1(\Omega))$.

There only remains to check that $\int_0^T
g(t,x;t/\e,x/e,\xi,\om)\:dt$ goes strongly to 0 in
$L^1_{\text{loc}}$ norm as $\e$ vanishes; this result follows
immediately from the same property for $g_n$ and the above
inequalities. Therefore, we skip its proof.

\vskip1mm

\begin{proof}[Proof of Lemma \ref{lem:finite_prop_speed}]

First, let us recall that
$$
\fe(t,x,\xi,\om)=f_0\left(\e Y\left(
\frac{t}{\e},\frac{x}{\e},\xi,\om\right),Y\left(
\frac{t}{\e},\frac{x}{\e},\xi,\om\right), \Xi\left(
\frac{t}{\e},\frac{x}{\e},\xi,\om\right),\om \right),
$$
and the Jacobian of the change of variables
$$
(x,\xi) \to \left(\e Y\left(
\frac{t}{\e},\frac{x}{\e},\xi,\om\right), \Xi\left(
\frac{t}{\e},\frac{x}{\e},\xi,\om\right),\om \right)
$$
is equal to 1.

On the other hand, since
$$
\frac 1 2 \left| \Xi(t,y,\xi,\om)\right|^2 + u\left(
Y(t,y,\xi,\om)\right)=\frac 1 2 |\xi|^2 + u(y,\om)
$$
we have
$$
\left| \Xi(t,y,\xi,\om)\right|\leq \sqrt{|\xi|^2 +2 u(y,\om)}\leq
\sqrt{|\xi|^2 +2 u_{\text{max}}}
$$
and
$$
\left| \e Y\left( \frac{t}{\e},\frac{x}{\e},\xi,\om\right)-
x\right|\leq t\sqrt{|\xi|^2 +2 u_{\text{max}}}.
$$
Thus
\begin{eqnarray*}
&&\int_{x\in B_R,\ \xi\in
B_{R'}}\left|\fe(t,x,\xi,\om)\right|\:dx\:d\xi\\&=&\int_{x\in
B_R,\ \xi\in B_{R'}}\left|f_0\left(\e Y\left(
\frac{t}{\e},\frac{x}{\e},\xi,\om\right),Y\left(
\frac{t}{\e},\frac{x}{\e},\xi,\om\right), \Xi\left(
\frac{t}{\e},\frac{x}{\e},\xi,\om\right),\om
\right)\right|\:dx\:d\xi\\
&\leq &\int_{x\in B_R,\ \xi\in B_{R'}}\sup_{y\in
\R^N}\left|f_0\left(\e Y\left(
\frac{t}{\e},\frac{x}{\e},\xi,\om\right),y, \Xi\left(
\frac{t}{\e},\frac{x}{\e},\xi,\om\right),\om
\right)\right|\:dx\:d\xi\\
&\leq & \int_{\R^N\times\R^N }\mathbf 1_{(\e Y\left(
\frac{t}{\e},\frac{x}{\e},\xi,\om\right),\Xi\left(
\frac{t}{\e},\frac{x}{\e},\xi,\om\right))\in K_{T,R,R'}}\sup_{y\in
\R^N}\left|f_0\left(\e Y\left(
\frac{t}{\e},\frac{x}{\e},\xi,\om\right),y, \Xi\left(
\frac{t}{\e},\frac{x}{\e},\xi,\om\right),\om
\right)\right|\:dx\:d\xi\\
&=&\int_{K_{T,R,R'}}\sup_{y\in \R^N}\left|f_0\left(x,y, \xi,\om
\right)\right|\:dx\:d\xi
\end{eqnarray*}

The proof of the other inequality goes along the same lines.

\end{proof}

\section{The integrable case}

In this section, we treat independently the periodic and the
stationary ergodic case. Indeed, some results of the periodic case
are no longer true in the stationary ergodic setting, and the
results which do remain valid are not proved with the same tools.

Let us make precise what we mean about ``integrable case'' : in
the periodic case, we take a function $u(y)$ which has the form
\be u(y)=\sum_{i=1}^N u_i(y_i), \label{hyp:int_per}\ee where each
function $u_i$ is periodic with period 1 ($1\leq i\leq N$). The
Hamiltonian $H(y,\xi)$ can be written
$$
H(y,\xi)=\frac 1 2 |\xi|^2 + u(y) =\sum_{i=1}^N H_i(y_i,\xi)
$$
where $H_i(y_i,\xi)=\frac 1 2 |\xi_i|^2 + u_i(y_i)$ ($1\leq i\leq
N$). And the Hamiltonian system \eqref{eq:systhamilt} becomes
\be
\left\{ \begin{array}{l} \dot{Y}_i = -\Xi_i,\\\dot{\Xi}_i
=u_i'(Y_i),\\
Y_i(t=0)=y_i, \quad \Xi_i(t=0)=\xi_i.
\end{array}\right.
\label{eq:systhamilt_onedim}\ee Thus it is enough to investigate
the behavior of each one-dimensional Hamiltonian system
\eqref{eq:systhamilt_onedim} individually, and for most
calculations, we can assume without loss of generality that $N=1$,
and we drop all indices $i$. However, for the calculation of the
projection $P$, a more thorough discussion will be needed, and we
will come back to the case when $N>1$ in the corresponding
paragraph.

In the stationary ergodic setting, expression \eqref{hyp:int_per}
can be transposed in the following way : assume that
$\Omega=\Pi_{i=1}^N \Omega_i$, where each $\Omega_i $ is a
probability space, and assume that for $1\leq i\leq N$, an ergodic
group transformation, denoted by $(\tau_{i,y})_{y\in\R}$, acts on
each $\Omega_i$.

Then for $\om=(\om_1,\cdots,\om_N)\in\Omega$, and
$y=(y_1,\cdots,y_N)\in\R^N$, we set
$\tau_y\om:=(\tau_{1,y_1}\om_1,\cdots, \tau_{N,y_N}\om_N)$. And we
assume that the function $u$ can be written
$$
u(y,\om)=\sum_{i=1}^N U_i\left(\tau_{i,y_i}\om_i \right),
$$
where $U_i\in L^{\infty}(\Omega_i)$ for all $1\leq i\leq N$. The
same remarks as in the periodic case can be made, and thus we will
only consider the case $N=1$; note that in the stationary ergodic
case, we are unable to compute the projection $P$ when $N>1$.

\subsection{Periodic setting}
The goal of this subsection is to give another proof of the
results of K. Hamdache and E. Frénod in \cite{HamdacheFrenod},
based on the study of the system \be \left\{ \begin{array}{l}
\dot{Y} = -\Xi,\\\dot{\Xi}
=u'(Y),\\
Y(t=0)=y, \quad \Xi=\xi, \quad y\in\R,\ \xi\in\R.
\end{array}\right.
\label{eq:systhamiltper}\ee

The Hamiltonian $H(y,\xi)=\frac 1 2 |\xi|^2 + u(y)$ is constant
along the trajectories of the system \eqref{eq:systhamiltper}, so
that
$$
\frac 1 2 |\Xi(t,y\xi)|^2 + u(Y(t,y,\xi))=H(y,\xi).
$$

We now fix $y,\xi\in\R^N$. Without any loss of generality, we
assume $y\in[-\frac 1 2,\frac 1 2)$, and we set $\mathcal
E:=H(y,\xi)$. The above equation describes the movement of a
single particle in a periodic potential $u$, with $0\leq u\leq
u_{\text{max}}$. It is well-known that there are two kinds of
behavior, depending on the value of the energy $\mathcal E$~: if
$\mathcal E<u_{\text{max}}$, the particle is ``trapped'' in a well
of potential around $y$, and $Y(t)$ remains bounded as
$t\to\infty$. In that case, the trajectories in the phase space
are closed curves. If $\mathcal E>u_{\text{max}}$, the trajectory
of the particle is unconstrained and $|Y(t)|\to \infty$ as
$t\to\infty$. We study more precisely these two cases and their
consequences on the expression of the projection $P$ in the
following; we refer for instance to \cite{arnold} for further
calculations and results about Hamiltonian dynamics and ordinary
differential equations in general.

\subsubsection{Expression of $\xs(y,\xi)$}

We begin with the case when $H(y,\xi)<u_{\text{max}}$. In that
case, $u(y)\leq H(y,\xi)<u_{\text{max}}$. By continuity of the
potential $u$, there exists $y_-<y$ and $y_+>y$ such that
$H(y,\xi)<u(y_{\pm})<u_{\text{max}}$, and the periodicity of $u$
allows us to choose $y_{\pm}$ such that $|y_+-y_-|<1$. Then
$y_-<Y(t,y,\xi)<y_+$ for all $t\geq 0$. Indeed, assume that there
exists $t>0$ such that $Y(t,y,\xi)\geq y_+>y=Y(t=0,y,\xi)$. Since
the trajectory $Y$ is continuous in time, there exists $0<t_0\leq
t$ such that $Y(t=t_0,y,\xi)=y_+$, which is absurd since
$$
H(Y(t=t_0,y,\xi), \Xi(t=t_0,y,\xi))=H(y,\xi)\geq
u(Y(t=t_0,y,\xi))>H(y,\xi).
$$
Thus $Y(t,y,\xi)$ is bounded. Since
$$
\xs(y,\xi)=\lim_{T\to\infty}\frac 1 T\int_0^T
\Xi(t,y,\xi)\:dt=-\lim_{T\to\infty}\frac 1 T\int_0^T
\dot{Y}(t,y,\xi)\:dt=\lim_{T\to\infty}\frac{y-Y(T,y,\xi)}{T}
$$
we deduce that $\xs(y,\xi)=0$ for all $y,\xi$ such that
$H(y,\xi)<u_{\text{max}}$.

\vskip2mm We now study the case $H(y,\xi)>u_{\text{max}}$. Since
$$
|\dot{Y}(t,y,\xi)|^2=2\left( H(y,\xi) - u(Y(t,y,\xi))\right)\geq 2
(H(y,\xi)-u_{\text{max}})>0
$$
$\dot{Y}$ does not vanish for $t\geq 0$. Consequently,
$$
\Xi(t,y,\xi)=-\dot{Y}(t,y,\xi)=\sgn(\xi) \sqrt{2\left( H(y,\xi) -
u(Y(t,y,\xi))\right)},
$$
and since $|Y(t,y,\xi)-y|\geq \sqrt{2\left(
H(y,\xi)-u_{\text{max}}\right)}t$, $|Y(t)|\to \infty$ as
$t\to\infty$. We immediately deduce that $\Xi(t,y,\xi)$ is
periodic in time: indeed, there exists $t_0>0$ such that
$$
Y(t_0,y,\xi)=y -\sgn (\xi).
$$
And $$\Xi(t=t_0,y,\xi)=\sgn(\xi) \sqrt{2\left( H(y,\xi) -
u(Y(t_0,y,\xi))\right)}=\sgn(\xi) \sqrt{2\left( H(y,\xi) -
u(y)\right)}=\xi=\Xi(t=0,y,\xi),$$ so that for $s\geq 0$,
\begin{gather*}
Y(t_0+s,y,\xi)=Y(s,y,\xi)-\sgn (\xi),\\
\Xi(t_0+s,y,\xi)=\Xi(s,y\xi),
\end{gather*}
and $\Xi$ is periodic with period $t_0$.

Consequently,
$$
\xs(y,\xi)=\lim_{T\to\infty}\frac 1 T \int_0^T\Xi(t,y,\xi)\:dt
=\frac 1 {t_0} \int_0^{t_0}\Xi(t,y,\xi)\:dt.
$$

But\begin{eqnarray*}
\int_0^{t_0}\Xi(t,y,\xi)\:dt&=&-\int_0^{t_0}\dot{Y}(t,y,\xi)\:dt\\
&=&-\left(Y(t_0,y,\xi)-y \right)\\
&=&\sgn (\xi).
\end{eqnarray*}
Thus we only have to compute $t_0$. With this aim in view, we use
the change of variables $s=Y(t)$, with Jacobian $ds=\dot{Y}dt$
(recall that $\dot{Y}(t,y,\xi)=-\sgn(\xi) \sqrt{2\left( H(y,\xi) -
u(Y(t,y,\xi))\right)}$\;), in the formula
\begin{eqnarray*}
t_0&=&\int_0^{t_0}dt\\
&=& \int_{Y(t=0)}^{Y(t_0)}\frac{1}{-\sgn(\xi) \sqrt{2\left(
H(y,\xi) - u(s)\right)}}ds\\
&=&-\sgn(\xi)\int_{y}^{y -\sgn (\xi)} \frac{1}{ \sqrt{2\left(
H(y,\xi) - u(s)\right)}}ds\\
&=&\int_0^{1}\frac{1}{ \sqrt{2\left( H(y,\xi) - u(s)\right)}}ds
\end{eqnarray*}

Eventually, we deduce
$$
\xs(y,\xi)=\sgn(\xi) \varphi(H(y,\xi)),
$$
where
$$
\varphi(\mathcal E)=\sqrt{2}\mathbf 1_{\mathcal
E>u_{\max}}\frac{1}{\mean{\frac{1}{\sqrt{\left( \mathcal E -
u(s)\right)}}}}
$$

We close this paragraph with a calculation which allows us to
express $\xs$ in terms of the homogenized Hamiltonian $\bar{H}$.
The result we will obtain will be justified in more abstract and
theoretical terms in the last subsection, using arguments similar
to those of the theory of Aubry-Mather.

First, let us recall the expression of the homogenized Hamiltonian
$\bar{H}$ (see \cite{LPV}) : we have
$$
H(y,\xi)=\frac  1 2 |\xi|^2 + u(y), \quad \text{with} \inf u=0,\
\sup u = u_{\text{max}},
$$
and thus
$$
\bar{H}(p)=u_{\text{max}} + \frac 1 2 \left\{\begin{array}{ll}
                                                0\quad&\text{if}\ p<
                                                \mean{\sqrt{2(u_{\text{max}}-u)}}\\
                                                \lambda\quad&\text{if}\
                                                |p|\geq
                                                \mean{\sqrt{2(u_{\text{max}}-u)}},\quad
                                                \text{where }|p|=
                                                \mean{\sqrt{2(u_{\text{max}}-u)+\lambda}}\end{array}\right.
$$

In other words, setting
$$
\theta:\begin{array}{ccc}[0,\infty)&\to&[0,\infty)\\
\lambda&\mapsto&\mean{\sqrt{2(u_{\text{max}}-u)+\lambda}}\end{array}
$$
we have
$$
\bar{H}(p)=u_{\text{max}} + \frac 1 2 \mathbf 1_{|p|\geq
\theta(0)}\theta^{-1}(|p|).
$$
Hence,
$$
\bar{H}'(p)=\sgn(p)\frac 1 2 \mathbf 1_{|p|\geq
\theta(0)}\frac{1}{\theta'\left(\theta^{-1}(|p|) \right)};
$$
and
\begin{gather*}
\theta'(\lambda)=\frac 1
2\mean{\frac{1}{\sqrt{2(u_{\text{max}}-u)+\lambda}}},\\
\theta^{-1}(|p|)=2\left(\bar{H}(p)-u_{\text{max}} \right)\quad
\forall |p|\geq \theta(0),\\
|p|> \theta(0)\iff \bar{H}(p)>u_{\max}\quad \forall p.
\end{gather*}
Gathering all the terms, we are led to
\begin{eqnarray*}
\bar{H}'(p)&=&\sgn(p) \sqrt{2}\mathbf 1_{
\bar{H}(p)>u_{\max}}\frac{1}{\mean{\frac{1}{\sqrt{\bar{H}(p)-u}}}}\\
&=&\sgn(p) \varphi\left( \bar{H}(p)\right) \end{eqnarray*}

Thus, the final expression is
$$
\xs(y,\xi)=\bar{H}'(p),
$$
where $p$ is such that
$$
\bar{H}(p)=H(y,\xi)\vee u_{max}, \quad \sgn(p)=\sgn(\xi).
$$
\label{sssec:xs_per}
\subsubsection{Expression of the projection $P$}

We also mention here how to find a general expression of the
projection $P$ in the special case $N=1$, and we explain how to
generalize this expression in some particular cases when $N>1$.
Recall that if $f=f(y,\xi)\in
L^1_{\text{loc}}(\R^N_y\times\R^N_{\xi})$ is periodic in $y$, then
$$
P(f)(y,\xi)=\lim_{T\to\infty}\frac 1 T \int_0^T
f(Y(t,y,\xi),\Xi(t,y,\xi))\:dt
$$
and the limit holds almost everywhere and in $L^1([0,1)\times\R^N,
m_c)$, with $dm_c(y,\xi)=\mathbf 1_{H(y,\xi)\leq c }\:dy \:d\xi$.

We begin with the case $H(y,\xi)>u_{\max}$. We have seen in the
previous paragraph that there exists $t_0>0$, which depends only
on $H(y,\xi)$ such that for all $t>0$, for all $k\in\N$
$$
Y(t+ k t_0,y,\xi)=Y(t,y,\xi ) -k\sgn (\xi) ,\quad  \Xi(t+ k
t_0,y,\xi)=\Xi(t,y,\xi ).
$$
Thus $f(Y(t), \Xi(t))$ is periodic in time with period $t_0$, and
$$
P(f)(y,\xi)=\frac{1}{t_0}\int_0^{t_0}f(Y(t,y,\xi),\Xi(t,y,\xi))\:dt.
$$
We use once again the change of variables $s=Y(t)$, so that
\begin{eqnarray*}
&&\int_0^{t_0}f(Y(t,y,\xi),\Xi(t,y,\xi))\:dt\\&=&\int_y^{y-\sgn(\xi)
} f(s,\sgn(\xi) \sqrt{2\left( H(y,\xi) - u(s)\right)})
\frac{1}{-\sgn(\xi) \sqrt{2\left(
H(y,\xi) - u(s)\right)}} \:ds\\
&=&\mean{f\left(\cdot,\sgn(\xi) \sqrt{2\left( H(y,\xi) -
u(\cdot)\right)}\;\right) \frac{1}{ \sqrt{2\left( H(y,\xi) -
u(\cdot)\right)}} }.
\end{eqnarray*}
And eventually, \be P(f)(y,\xi)=\bar{f}(\sgn(\xi),
H(y,\xi))\label{eq:proj_per_1} \ee with
$$
\bar{f}(\eta,\mathcal E):=\frac{\mean{f\left(\cdot,\eta
\sqrt{2\left( \mathcal E - u(\cdot)\right)}\;\right) \frac{1}{
\sqrt{\left(\mathcal E - u(\cdot)\right)}} }}{\mean{\frac{1}{
\sqrt{\left( \mathcal E - u\right)}}}}\quad \eta=\pm 1,\ \mathcal
E>u_{\text{max}}
$$

We now focus on the case $0<\mathcal E<u_{\text{max}}$. In order
to simplify the analysis we assume that $\mathcal E\notin \{u(y)\;
; u\ \text{has a local extremum at }y\}$ (this set is finite or
countable), and that
$$
\forall y\in\R,\quad u'(y)=0\Rightarrow u \ \text{has a local
extremum at }y.
$$

In that case, it can be easily proved that $Y(t,y,\xi)$ is
periodic in $t$; this follows directly from the fact that the
trajectory in the phase space is closed (see \cite{arnold}).
Indeed, pushing a little further the analysis of the previous
paragraph, we construct $z_{\pm}$ such that
\begin{gather*}
|z_+-z_-|<2\pi,\quad z_-<z_+,\\
u(z_{\pm})=\mathcal E,\\
z_-\leq y\leq z_+,\\
u(z)<\mathcal E\quad \forall z\in(z_-,z_+).
\end{gather*}
Then the particle starting from $y$ with initial speed $-\xi$
reaches either $z_+$ or $z_-$ in finite time; the speed of the
particle is 0 at that moment since
$$
|\dot{Y}|^2=2(\mathcal E-u(Y)),
$$
but its acceleration is $-u'(z_{\pm})\neq 0$, so the particle
turns around and goes back in the reverse direction. It then
reaches the other extremity of the interval $(z_-,z_+)$ in finite
time, and the same phenomena occurs. Hence after a finite time
$t_0$, the particle is back at its starting point $y$ with the
same speed $-\xi$. Consequently, the movement of the particle is
periodic in time with period $t_0$. Thus, we have
$$
P(f)(y,\xi)=\frac{1}{t_0}\int_{t_1}^{t_1+t_0}f(Y(t,y,\xi),\Xi(t,y,\xi))\:dt,
$$
where $t_1\geq 0$ is arbitrary. It is convenient to choose for
$t_1$ the first time when the particle hits $z_-$. In that case,
it is easily seen that $t_0$ is twice the time it takes to the
particle to go from $z_- $ to $z_+$, so that
$$
\frac{t_0}{2}=\int_{t_1}^{t_1+t_0/2}dt=\int_{z_-}^{z_+}\frac{1}{
\sqrt{2\left(\mathcal E - u(s)\right)}}ds= \mean{\mathbf
1_{u<\mathcal E}\frac{1}{ \sqrt{2\left( \mathcal E - u\right)}}}
$$
and
\begin{gather*}
\int_{t_1}^{t_1+\frac{t_0}{2}}f(Y(t,y,\xi),\Xi(t,y,\xi))\:dt=
\mean{\mathbf 1_{u<\mathcal E}f(s,-\frac{1}{ \sqrt{2\left(
\mathcal E -
u\right)}})\frac{1}{ \sqrt{2\left( \mathcal E - u\right)}}},\\
\int_{t_1+\frac{t_0}{2}}^{t_1+t_0}f(Y(t,y,\xi),\Xi(t,y,\xi))\:dt =
\mean{\mathbf 1_{u<\mathcal E}f(s,\frac{1}{ \sqrt{2\left( \mathcal
E - u\right)}})\frac{1}{ \sqrt{2\left( \mathcal E - u\right)}}}.
\end{gather*}
Gathering all the terms, we are led to \be
P(f)(y,\xi)=\frac{\mean{\mathbf 1_{u<\mathcal
E}\left[f\left(\cdot,\frac{1}{ \sqrt{2\left( \mathcal E -
u\right)}}\right) +f\left(\cdot,-\frac{1}{ \sqrt{2\left( \mathcal
E - u\right)}}\right)\right]\frac{1}{ \sqrt{\left( \mathcal E -
u\right)}}}}{2\mean{\mathbf 1_{u<\mathcal E}\frac{1}{ \sqrt{\left(
\mathcal E - u\right)}}}} \label{eq:proj_per_2}\ee

Expressions \eqref{eq:proj_per_1} and \eqref{eq:proj_per_2} are
compatible with the ones in \cite{HamdacheFrenod}.

\vskip2mm

Let us now come back to the case when $N>1$, and take a function
$\varphi(y,\xi)=\varphi_1(y_1,\xi_1)\cdots \varphi_N(y_N,\xi_N)$,
where each $\varphi_i$ is periodic with period 1. We want to
compute the limit
$$
\frac 1 T \int^T_0
\varphi_1(Y_1(t,y_1,\xi_1),\Xi_1(t,y_1,\xi_1))\cdots
\varphi_N(Y_N(t,y_N,\xi_N),\Xi_N(t,y_N,\xi_N))\:dt.
$$
In general, knowing the behavior of each trajectory $(Y_i,\Xi_i)$
independently is not enough to compute such a product. However,
here, we recall that each function
$\varphi_i(Y_i(t,y_i,\xi_i),\Xi_i(t,y_i,\xi_i))$ ($1\leq i\leq N$)
is periodic in time. The period depends only on $H_i(y_i,\xi_i)$
and on the function $u_i$. More precisely, setting
$$
T_i(\mathcal E):=\sqrt 2\int_0^1\mathbf 1_{u_i(z)<\mathcal
E}\frac{1}{\sqrt{\mathcal E-u_i(z)}}\:dz\quad \forall \mathcal
E>0, \mathcal E\neq u_{\text{max}},
$$
$\varphi_i(Y_i(t,y_i,\xi_i),\Xi_i(t,y_i,\xi_i))$ is periodic in
time with period $T_i(H_i(y_i,\xi_i))$.

We can thus use the following result :
\begin{lem}
Let $f_1,\cdots,f_N\in L^{\infty}(\R)$ such that $f_i$ is periodic
with period $\theta_i$, and set $ \mean{f_i}=\frac 1 {\theta_i }
\int_0^{\theta_i}f_i.$

Assume that \be \frac{k_1}{ \theta_1} + \cdots + \frac{k_N}{
\theta_N}\neq 0\quad \forall (k_1,\cdots,k_N)\in\Z^N\setminus\{
0\}.\label{hyp:noncommens}\ee

Then as $T\to\infty$,
$$
\frac 1 T \int_0^T f_1(t)\cdots f_N(t)\:dt \to
\mean{f_1}\cdots\mean{f_N}.
$$
\end{lem}
\begin{proof}[Sketch of proof]
By density, it is enough to prove the lemma for $f_1,\cdot, f_N\in
\mathcal C^{\infty}(\R)$. Write $f_i$ as a Fourier series (the
series converges thanks to the regularity assumption), and use the
fact that for all $\alpha\neq 0$,
$$
\frac 1 T \int_0^T e^{i\alpha t}\:dt\to 0\quad \text{as
}T\to\infty.
$$

\end{proof}

In the present setting, we deduce the following result :

\begin{prop}

Let $\varphi:(y,\xi)\mapsto\varphi_1(y_1,\xi_1)\cdots
\varphi_N(y_N,\xi_N)$, where $\varphi_i\in
L^{\infty}_{\text{per}}(\R_y\times\R^{\xi})$.

Let $(y,\xi)\in[0,1)^N\times \R^N$, and let
$\theta_i=\theta_i(y_i,\xi_i)=T_i(H_i(y_i,\xi_i))$ for $1\leq
i\leq N$. Assume that $(\theta_1,\cdots,\theta_N)$ satisfy
condition \eqref{hyp:noncommens}. Then \be
P(\varphi)(y,\xi)=P_1(\varphi_1)(y_1,\xi_1)\cdots
P_N(\varphi_N)(y_N,\xi_N)\label{eq:proj_per_N} \ee where each
$P_i$ is the projection in dimension 1 with potential $u_i$, given
by expressions \eqref{eq:proj_per_1} and \eqref{eq:proj_per_2}.

\end{prop}

In particular, when the set
$$
\{(y,\xi)\in[0,1]^N\times\R^N;
\;(\theta_1(y_1,\xi_1),\cdots,\theta_N(y_N,\xi_N))\ \text{satisfy
condition}\ \eqref{hyp:noncommens} \}
$$
has zero Lebesgue measure, equality \eqref{eq:proj_per_N} holds
almost everywhere. It can then be generalized to arbitrary
functions $\varphi\in L^{\infty}_{\text{per}}(\R^N\times\R^N)$
(always by linearity and density). The correct expression of the
projection $P$ is then \be P=P_1\circ P_2\circ\cdots\circ P_N,
\label{eq:proj}\ee where each projection $P_i$ acts on the
variables $(y_i,\xi_i)$ only. Notice that all projections $P_i$
thus commute with one another; hence the order in which they are
taken is unimportant.

We wish to emphasize that on the open set $\{(y,\xi)\in\R^{2N},
\forall i\in\{1,\cdots,N\}\ H_i(y_i,\xi_i)>\max u_{i}\}$, the
expression \eqref{eq:proj} is true. Indeed, the function $T_i$ is
strictly decreasing on $(\max u_{i}, +\infty)$, and thus the set
$$
\{(\mathcal E_1, \cdots,\mathcal E_N)\in\R^N;\ \mathcal E_i>\max
u_{i}\ \text{and } \exists k\in\Z^N\setminus \{0\},\ k_1/
T_1(\mathcal E_1) + \cdots + k_N/ T_N(\mathcal E_N)=0  \}
$$
is countable. As a consequence, the set $$\{(y,\xi)\in\R^{2N},
H_i(y_i,\xi_i)>\max u_{i}\ \forall i\ \text{and}\
(\theta_1(y_1,\xi_1),\cdots,\theta_N(y_N,\xi_N))\ \text{satisfy
condition}\ \eqref{hyp:noncommens}  \}$$ has zero Lebesgue
measure.

\vskip1mm

However, let us mention here that in general, condition
\eqref{hyp:noncommens} cannot be relaxed : indeed, assume for
instance that $u_i=u_j:=u$ for $i\neq j$ and assume that the
function $u$ is such that
$$
\exists y_0>0,\quad u(y)=y^2\ \text{for}\ |y|<y_0,
$$
and $u(y)>y_0^2$ if $y\in [-\frac 1 2, \frac 1 2]\setminus [-y_0,
y_0]$.

Then if $|\mathcal E|\leq \sqrt{y_0}$, we have
$$
T(\mathcal E)=\int_{-\sqrt{\mathcal E}}^{\sqrt{\mathcal
E}}\frac{1}{\sqrt{\mathcal E-y^2}}\ dy = 2 \int_0^1\frac 1
{\sqrt{1-z^2}}\:dz=:T_0
$$
Thus, if $H_i(y_i,\xi_i)\leq\sqrt{y_0}$, then
$(Y_i,\Xi_i)(t,y_i,\xi_i)$ is periodic with period $T_0$. Notice
that $T_0$ does not depend on the energy $H_i(y_i,\xi_i)$

In that case, the function $\varphi(Y(t),\Xi(t))$ is also periodic
with period $T_0$. Thus we have to compute the limit of
$$
\frac 1 T \int_0^T f_1(t)\cdots f_N(t)\: dt
$$
as $T\to\infty$, where the $f_i$ are arbitrary functions with
period $T_0$. It is then easily proved that \be \frac 1 T \int_0^T
f_1(t)\cdots f_N(t)\: dt\to \sum_{\begin{array}{c}k\in\Z^N,\\k_1 +
\cdots +k_N=0\end{array}}a_{1,k_1}\cdots a_{N,k_N}
\label{neg_result}\ee where
$$
a_{j,l}=\frac 1 {T_0}\int_0^{T_0}f_j(t) e^{-\frac{2il\pi
t}{T_0}}\:dt,\quad 1\leq j\leq N,\ l\in\Z.
$$
In general, the right-hand side of \eqref{neg_result} differs from
$a_{1,0}\cdots a_{N,0}$, and thus
$$
P\neq P_1\circ\cdots \circ P_N
$$
for $(y,\xi)$ in a neighbourhood of the origin.

\subsection{Stationary ergodic setting}

In the stationary ergodic setting, some of the expressions or
properties above are no longer true. The most significant
difference occurs when the energy $H(y,\xi)<u_{\text{max}}$;
indeed, in that case the particle is not necessarily trapped,
depending on the profile of the potential $u$. Hence, in the rest
of the subsection, we focus on the case $H(y,\xi)>u_{\text{max}}$.
In that case, the movement of the particle is unbounded and has
many similarities with the periodic case. In particular, the
particle sees ``all the potential'' during its evolution, and this
will be fundamental in the use of the ergodic theorem.

\subsubsection{Expression of $\xs(y,\xi, \om)$}

This paragraph is devoted to the proof of proposition
\ref{prop:xs_ergo} in the stationary ergodic setting. We refer for
instance to \cite{LionsSougCorrectors} for conditions on the
existence of correctors for all $P\in\R$ in the case of a general
coercive hamiltonian. In the present case, there exist correctors
if
$$
\{ y\in\R; \ u(y,\om)=u_{\text{max}}\}\neq \emptyset
$$
a.s. in $\om\in\Omega$.

\begin{rem}
We wish to point out that the expressions in the periodic and in
the stationary ergodic case when $H(y,\xi,\om)>u_{\text{max}}$ are
exactly the same (compare proposition \ref{prop:xs_ergo} and the
end of paragraph \ref{sssec:xs_per}). This expression, and more
precisely, the equality $\xs=\bar{H}'(P)$ for some $P$, is in fact
strongly linked to Aubry-Mather theory. Indeed,
$$
\xs(y,\xi,\om)=\lim_{T\to\infty} \frac 1 T \int_0^T
\Xi(t,y,\xi,\om)\: dt = -\lim_{T\to
\infty}\frac{Y(T,y,\xi,\om)-y}{T},
$$
and $\xs(y,\xi,\om)$ is thus (up to a multiplication by $-1$) the
rotation number associated to the Hamiltonian flow starting at
$(y,\xi)$. The interested reader should compare our proposition
\ref{prop:xs_ergo} to lemma 2.8 in \cite{WEAubryMather} or theorem
4.1 in \cite{EvansGomes1}, and our proof to the ones in these
articles. We refer to \cite{WEAubryMather,EvansGomes1} for further
references to Aubry-Mather theory and its applications to partial
differential equations.
\end{rem}

\begin{proof}[Proof of proposition \ref{prop:xs_ergo}]

In all the proof, we fix $y,\xi,\om$ such that
$H(y,\xi,\om)>u_{\text{max}}$, and we set $P=P(y,\xi,\om)$. Let
$Q\in\R$ be arbitrary, and let $v$ be a corrector, i.e.
$$
H(y, Q +\nabla_y v(y,\om),\om)=\bar{H}(Q).
$$

Then according to the theory of viscosity solutions, for all
$(y,\xi,\om)\in\R\times\R\times\Omega$, for all $T>0$, \be
v(y,\om) \leq v(Y(T,y,\xi,\om)) + \int_0^T
L(Y(t,y,\xi,\om),\Xi(t,y,\xi,\om),\om) \;dt +
Q\left[Y(T,y,\xi,\om) -y \right] + \bar{H}(Q) T.
\label{in:viscosity}\ee

Hence \be \frac 1 T \int_0^T
L(Y(t,y,\xi,\om),\Xi(t,y,\xi,\om),\om) \;dt\geq \frac{v(y)
-v(Y(T,y,\xi,\om))}{T} - Q \frac{Y(T,y,\xi,\om) -y }{T} -
\bar{H}(Q) T \label{in:ergo_L}\ee and
$$
\frac 1 2 |\dot{Y}(t,y,\xi,\om)|= H(y,\xi,\om) - u \geq
H(y,\xi,\om) - u_{\text{max}}>0\quad \forall t>0.
$$
Thus there exist constants $\alpha,\beta>0$ depending only on
$H(y,\xi,\om)$ and $u_{\text{max}}$, such that
$$
0<\alpha\leq \left|\frac{Y(T,y,\xi,\om) -y }{T} \right|\leq
\beta\quad \forall T>0.
$$
Consequently, $Y(T)\to \infty$ as $T\to\infty$ and
$$
\frac{v(y) -v(Y(T,y,\xi,\om))}{T} =\frac{v(y)
-v(Y(T,y,\xi,\om))}{Y(T,y,\xi,\om) -y }\;\frac{Y(T,y,\xi,\om) -y
}{T}\to 0\quad \text{as}\ T\to\infty.
$$
(remember \eqref{hyp:correctsublin}).

On the other hand, as $T\to \infty$,
$$
\frac{Y(T,y,\xi,\om) -y }{T}=-\frac 1 T \int_0^T
\Xi(t,y,\xi,\om)\:dt\to -\xs(y,\xi,\om).
$$
Hence, passing to the limit in \eqref{in:ergo_L}, we derive \be
P(L)(y,\xi,\om)\geq Q\xs(y,\xi,\om) - \bar{H}(Q)\quad \forall
Q\in\R. \ee

Thus \be P(L)(y,\xi,\om)\geq \bar{L}\left( \xs(y,\xi,\om)\right).
\label{in:proj_L}\ee

\vskip1mm

In order to prove the proposition, we have to find a special
$Q_0\in\R$ such that equality holds in \eqref{in:viscosity}. This
will entail that
$$
P(L)(y,\xi,\om)=\sup_{Q\in\R} \left(Q\xs(y,\xi,\om) -
\bar{H}(Q)\right)= \bar{L}\left( \xs(y,\xi,\om)\right),
$$
and the $\sup$ is obtained for $\xs(y,\xi,\om)=\bar H'(Q_0)$.

Let us thus prove that with $Q=P=P(y,\xi,\om)$, equality holds in
\eqref{in:viscosity}.

First, notice that
$$
v(y,\om):=\sgn(P)\int_0^y\sqrt{2(\bar{H}(P)-u(z,\om))}\:dz - Py
$$
is a viscosity solution of
$$
H(y, P + \nabla_y v,\om)=\bar{H}(P),
$$
and as $y\to \infty$
$$
\frac{1}{y}\int_0^y\sqrt{2(\bar{H}(P)-u(z,\om))}\:dz\to
E\left[\sqrt{2(\bar{H}(P)-U)} \right].
$$

By definition of $\bar{H}$,
$$
E\left[\sqrt{2(\bar{H}(P)-U)} \right]=|P|;
$$
consequently,
$$
\frac{1}{1 + |y|
}\left(\sgn(P)\int_0^y\sqrt{2(\bar{H}(P)-u(z,\om))}\:dz -
Py\right)\to 0
$$
as $y\to\infty$, a.s. in $\om$. Thus $v$ satisfies
\eqref{hyp:correctsublin}, and $v\in L^{\infty}(\Omega; \mathcal
C^1(\R^N))$. Thus the method of characteristics, for instance, can
be used to prove that equality holds in \eqref{in:viscosity}, with
$(y,\xi)$ replaced by any couple $(y',\xi')=(y', P + \nabla_y
v(y',\om))$. We have to prove that we can take
$(y',\xi')=(y,\xi)$. First, notice that
$$
\nabla_y v(y',\om)=\sgn(P) \sqrt{2(\bar{H}(P)-u(y',\om))} - P,
$$
and thus $\sgn\left(P +\nabla_y
v(y',\om)\right)=\sgn(\xi')=\sgn(P)=\sgn(\xi)$. Hence, take
$y'=y$. Then $|\xi|^2=|\xi'|^2$ because
$H(y,\xi)=\bar{H}(P)=H(y,\xi')$ by definition of $\xi'$. Thus
$\xi=\xi'$, and equality holds in \eqref{in:viscosity}.

\end{proof}

\subsection{Expression of the projection $P$}

The same method as in the periodic case can be used in order to
find the expression of the projection $P$ when
$H(y,\xi,\om)=:\mathcal E>u_{\text{max}}$; indeed, in that case,
remember that
$$
P(f)(y,\xi,\om)=\lim_{T\to\infty} \frac 1 T \int_0^T f\left(
Y(t,y,\xi,\om), \Xi(t,y,\xi,\om),\om\right)\:dt
$$
and we can use the change of variables
$$
dt=\frac{1}{\dot{Y}}\:dY=\frac{1}{-\sgn(\xi) \sqrt{2(\mathcal E -
u(Y,\om))}}\:dY
$$
in order to obtain
\begin{eqnarray*}
&&\int_0^T f\left( Y(t,y,\xi,\om),
\Xi(t,y,\xi,\om),\om\right)\:dt\\&=&\int_y^{Y(T,y,\xi,\om)}f\left(
z, \sgn(\xi) \sqrt{2(\mathcal E -
u(z,\om))},\om\right)\frac{1}{-\sgn(\xi) \sqrt{2(\mathcal E -
u(z,\om))}}\:dz. \end{eqnarray*} Since the group transformation
$(\tau_x)$ is ergodic, and $Y(T)\to\infty$ as $T\to\infty$, for
all $\mathcal E>u_{\text{max}}$,
\begin{multline*}
\frac{1}{Y(T)-y}\int_y^{Y(T,y,\xi,\om)}f\left( z, \sgn(\xi)
\sqrt{2(\mathcal E - u(z,\om))},\om\right)\frac{1}{-\sgn(\xi)
\sqrt{2(\mathcal E  - u(z,\om))}}\:dz\to\\\to E\left[F\left(
\sgn(\xi) \sqrt{2(\mathcal E  -
u(0,\om))},\om\right)\frac{1}{-\sgn(\xi) \sqrt{2(\mathcal E  -
u(0,\om))}} \right]
\end{multline*}

Thus, we obtain
$$
P(f)(y,\xi,\om)= \xs(y,\xi,\om)\bar{f}(sgn(\xi),H(y,\xi,\om)),
$$
where $$\bar{f}(\eta,\mathcal E)=E\left[F\left( \eta
\sqrt{2(\mathcal E  - u(0,\om))},\om\right)\frac{1}{\eta
\sqrt{2(\mathcal E  - u(0,\om))}} \right].$$

\bibliography{../articles,../books}

\end{document}